\documentclass[12pt]{article}
\usepackage{amsmath}
\usepackage{amsfonts}
\usepackage{amssymb}
\usepackage{theorem} 

\title{Polish $G$-spaces, the generalized model theory and complexity}
\author{A. Ivanov
, B. Majcher-Iwanow
 }
\date{ } 

\newtheorem{thm}{Theorem}[section]
\newtheorem{lem}[thm]{Lemma}
\newtheorem{definicja}[thm]{Definition}
\newtheorem{cor}[thm]{Corollary}
\newtheorem{prop}[thm]{Proposition}
\newtheorem{remark}[thm]{Remark}  

\begin{document}
\maketitle
\topskip 20pt

\begin{quote}
{\bf Abstract.} 
Given  Polish space ${\bf Y}$ and a continuous language $L$ 
we study the corresponding logic $\mathsf{Iso}({\bf Y})$-space ${\bf Y}_L$. 
We build a framework of generalized model theory towards analysis of 
Borel/algorithmic complexity of subsets of 
${\bf Y}^k_L\times  (\mathsf{Iso} ({\bf Y}))^l$. 

\bigskip

{\em 2010 Mathematics Subject Classification:} 03E15, 03C07, 03C57

{\em Keywords:}  Polish G-spaces, Continuous logic, Generalized model theory.
\end{quote}

\section{Introduction}

Let $({\bf Y},d)$ be a Polish space and $\mathsf{Iso}({\bf Y})$ 
be the corresponding isometry group 
endowed with the pointwise convergence topology. 
Then $\mathsf{Iso} ({\bf Y})$ is a Polish group. \parskip0pt

For any countable continuous signature $L$ the set 
${\bf Y}_L$ of all continuous metric $L$-structures on 
$({\bf Y},d)$ can be considered as a Polish 
$\mathsf{Iso}({\bf Y})$-space. 
We call this action {\em logic}. 
It is known that the logic action is universal 
for Borel reducibility of orbit equivalence relations 
of Polish $G$-spaces with closed $G\le \mathsf{Iso}({\bf Y})$ 
\cite{CL}, \cite{IMI-APAL}.  
Moreover by a result of 
J.Melleray in \cite{melleray} 
every Polish group $G$ can be 
realised as the automorphism 
group of a continuous metric structure on
an appropriate Polish space  $({\bf Y},d)$
and the structure is approximately ultrahomogeneous.  

On the other hand typical notions naturally arising for logic actions  
can be applied in the general case of a Polish 
$G$-space ${\bf X}$ with $G$ as above. 
If we consider $G$ together with a family of grey 
subgroups, then distinguishing an appropriate family 
$\mathcal{B}$ of grey subsets of ${\bf X}$ we arrive at 
the situation very similar to the logic space ${\bf Y}_L$ 
(see \cite{becker}, \cite{becker2}, \cite{MI05}, \cite{MI08}  
for non-Archimedian $G$ and \cite{IMI-APAL} for the general case). 
For example we can treat elements of $\mathcal{B}$ as 
continuous formulas. 
Then many theorems of traditional model theory can be 
generalized to topological statements concerning 
spaces with {\bf nice topologies} (i.e. defined by $\mathcal{B}$).
This approach is called the {\bf generalized model theory}, \cite{becker2}. 

\bigskip

The aim of this paper is to demonstrate that the tools 
of the generalized model theory nicely work for some  
other aspects of logic actions.  
Our basic concern is as follows. 
Viewing the logic space ${\bf Y}_L$ 
as a Polish space and using the recipe of generalized model theory 
we distinguish some subsets of 
${\bf Y}^k_L\times  (\mathsf{Iso} ({\bf Y}))^l$ and then  
study Borel/algorithmic complexity of them. 

We usually fix a countable dense subset $S_{{\bf Y}}$ of ${\bf Y}$ 
and study subsets of  ${\bf Y}_L$ which are invariant with respect to 
isometries stabilizing $S_{{\bf Y}}$ setwise.  
This is the approach where a structure on ${\bf Y}$ 
(say $M$) is considered together with its {\bf presentation over} $S_{{\bf Y}}$, 
i.e. together with the set 
$$ 
\mathsf{Diag} (M,S_{{\bf Y}})=\{ (\phi ,q): M\models \phi <q, \mbox{ where } q\in [0,1]\cap \mathbb{Q} 
\mbox{ and } \phi \mbox{ is a continuous } 
$$ 
$$ 
\mbox{ sentence with parameters from }S_{{\bf Y}} \}.     
$$  
The best example of this situation is the logic space $\mathfrak{U}_L$ over 
the  bounded Urysohn metric space $\mathfrak{U}$ 
\footnote{in this paper we do not use the Urysohn space $\mathbb{U}$;   
we only use the ball of it of diameter 1 and denote it by $\mathfrak{U}$} 
where distinguishing the countable counterpart $\mathbb{Q}\mathfrak{U}$ of $\mathfrak{U}$ 
(see Section 3 of \cite{IMI-APAL}) 
we study $\mathsf{Iso} (\mathbb{Q}\mathfrak{U})$-invariant subsets of $\mathfrak{U}_L$. 

We demonstrate in Section 2 that this setting appears as a part of 
generalized model theory studied in \cite{becker2} 
and \cite{IMI-APAL}. 
Moreover our methods also work in the Hilbert space case and 
in the measure algebras case. 

In Section 2.4 we consider 
other examples suggested by  
topological notions involving 
nice topologies.  
Since they correspond to natural logic constructions 
we view this material as a further development of generalized model theory. 

In Section 3 we show that our approach gives a framework  
to computable members of ${\bf Y}_L$ and their computable indexations.  
We will see that it supports standard approaches 
both to computable model theory and to effective metric spaces. 
Moreover it is suited to the general setting of computable Polish group actions 
presented in the recent paper \cite{MelMon}. 

In Section 4 we examine our approach in the cases of 
separable categoricity and ultrahomogeneity. 
In particular in Section 4.1 we find a Borel subset $\mathcal{SC}$ of 
${\bf Y}_L$ which is $\mathsf{Iso} (S_{{\bf Y}})$-invariant and 
can be viewed as the set of all presentations over $S_{{\bf Y}}$
of separably categorical structures on ${\bf Y}$.  
Moreover in Section 4.2 we study complexity of the index set 
of computable members of $\mathcal{SC}$. 

The paper is self-contained. 
We address it to logicians and do not assume any special background. 
We believe that the ideas presented in it can be 
helpful in model theory, descriptive set theory and computability theory.

\section{Logic space of continuous structures} 

Our paper belongs to a field of modern logic 
which can be situated between 
Invariant Descriptive Set Theory (\cite{bk}, \cite{gao-book}) 
and  Continuous Model Theory (\cite{BYBHU}, \cite{BNT}, \cite{BYU}), 
see also \cite{BYBM}, \cite{Chan} - \cite{CL}, \cite{RZ}. 
The definitions below basically correspond to these sources. 

\subsection{Polish group actions.}

A {\bf Polish space (group)} is a separable, completely
metrizable topological space (group).
Sometimes we extend the corresponding metric to 
tuples by 
$$
d((x_1 ,...,x_n ), (y_1 ,...,y_n ))= \mathsf{max} (d(x_1 ,y_1 ),...,d(x_n ,y_n )). 
$$ 
If a Polish group $G$ continuously acts on a Polish space ${\bf X}$,
then we say that ${\bf X}$ is a {\bf Polish $G$-space}. 
We  say that a subset of ${\bf X}$ is {\bf invariant} if
it is $G$-invariant. \parskip0pt

Let $({\bf Y},d)$ be a Polish space and $\mathsf{Iso}({\bf Y})$ 
be the corresponding isometry group 
endowed with the pointwise convergence toplogy. 
Then $\mathsf{Iso} ({\bf Y})$ is a Polish group. 
A compatible left-invariant metric can be obtained as follows: 
fix a countable dense set $S=\{ s_i : i\in \{ 1,2,...\} \}$ 
and then define for two isometries $\alpha$ 
and $\beta$ of ${\bf Y}$ 
$$ 
\rho_{S} (\alpha ,\beta )= \sum_{i=1}^{\infty} 2^{-i} \mathsf{min}(1, d(\alpha (s_i ),\beta (s_i ))) .
$$ 
 
We will study closed subgroups of $\mathsf{Iso}({\bf Y})$.  
We fix a dense countable set $\Upsilon \subset \mathsf{Iso}({\bf Y})$. 
In any closed subgroups of $\mathsf{Iso}({\bf Y},d)$ we distinguish 
the base consisting of all sets of the form 
$N_{\sigma ,q} = \{ \alpha : \rho_S (\alpha, \sigma )<q \}$, $\sigma \in \Upsilon$ 
and $q\in \mathbb{Q}$. 
We may assume that $\Upsilon$ is a subgroup of $\mathsf{Iso} ({\bf Y})$. 
To get this it is enough to replace  $\Upsilon$ by $G_0 = \langle  \Upsilon\rangle$.

\subsection{Continuous structures.} 

We now fix a countable continuous signature 
$$
L=\{ d,R_1 ,...,R_k ,..., F_1 ,..., F_l ,...\} . 
$$ 
Let us recall that a {\bf metric $L$-structure} 
is a complete metric space $(M,d)$ with $d$ bounded by 1, 
along with a family of uniformly continuous operations on $M$ 
and a family of predicates $R_i$, i.e. uniformly continuous maps 
from appropriate $M^{k_i}$ to $[0,1]$.   
It is usually assumed that $L$ assigns to each predicate symbol $R_i$ 
a continuity modulus $\gamma_i: [0,1] \rightarrow [0,1]$ so that 
any metric structure $M$ of the signature $L$ satisfies the property that if  
$d(x_j ,x'_j ) <\gamma_i (\varepsilon )$ with $1\le j\le k_i$,  
then the inequality  
$$ 
|R_i (x_1 ,...,x_j ,...,x_{k_i}) - R_i (x_1 ,...,x'_j ,...,x_{k_i})| < \varepsilon . 
$$ 
holds for the corresponding predicate of $M$. 
It happens very often that $\gamma_i$ coincides with $\mathsf{id}$. 
In this case we do not mention the appropriate modulus. 
We also fix continuity moduli for functional symbols. 

Note that each countable structure can be considered 
as a complete metric structure with the discrete $\{ 0,1\}$-metric.  

Atomic formulas are the expressions of the form $R_i (t_1 ,...,t_r )$, 
$d(t_1 ,t_2 )$, where $t_i$ are simply classical terms 
(built from functional $L$-symbols). 
We define {\bf formulas} to be expressions built from 
0,1 and atomic formulas by applications of the following functions: 
$$ 
x/2  \mbox{ , } x\dot- y= \mathsf{max} (x-y,0) \mbox{ , } \mathsf{min}(x ,y )  \mbox{ , } \mathsf{max}(x ,y )
\mbox{ , } |x-y| \mbox{ , } 
$$ 
$$ 
\neg (x) =1-x \mbox{ , } x\dot+ y= \mathsf{min}(x+y, 1) \mbox{ , } 
\mathsf{sup}_x \mbox{ and } \mathsf{inf}_x . 
$$ 
{\bf Statements} concerning metric structures are usually 
formulated in the form 
$$
\phi = 0, 
$$ 
where $\phi$ is a formula.  
Sometimes statements are called {\bf conditions}; we will use both names. 
A {\bf theory} is a set of statements without free variables 
(here $\mathsf{sup}_x$ and $\mathsf{inf}_x$ play the role of quantifiers). 

We often extend the set of formulas by the application 
of {\bf truncated products} by positive rational numbers.   
This means that when $q\cdot x$ is greater than $1$, the 
truncated product of $q$ and $x$ is $1$. 
Since the context is always clear, we preserve the same notation $q\cdot x$. 
The continuous logic after this extension does not differ from 
the basic case. 
   
It is worth noting that the choice of the set of connectives 
guarantees that for any continuous relational structure $M$, 
any formula $\phi$ is a $\gamma$-uniform 
continuous function from the 
appropriate power of $M$ to $[0,1]$, where 
$\gamma (\varepsilon )$ is of the form  
$$ 
\frac{1}{n} \cdot  
\mbox{min} \{ \gamma ' (\varepsilon ) : \gamma ' 
\mbox{ is a continuity modulus of an } 
L\mbox{-symbol appearing in the formula} \} , 
$$ 
$$ 
\mbox{ where the number $n$ only depends on the complexity of  }\phi . 
$$ 
This follows from the fact that when $\phi_1$ and $\phi_2$ have 
continuity moduli $\gamma_1$ and $\gamma_2$ respectively, 
then the formula $f(\phi_1 ,\phi_2 )$ obtained by applying 
a binary connective $f$, has a continuity modulus 
of the form $\mathsf{min} (\gamma_1 (\frac{1}{2} x) , \gamma_2 (\frac{1}{2} x) )$. 

It is observed in Appendix A of \cite{BYU} that instead of continuity moduli 
one can consider {\bf inverse continuity moduli}. 
Slightly modifying that place in \cite{BYU} we define it as follows. 

\begin{definicja} \label{inverse} 
A continuous monotone function $\delta :[0,1]\rightarrow [0,1]$ with $\delta(0) =0$ 
is an inverse  continuity modulus of a map $F(\bar{x}) :  {\bf X}^n \rightarrow [0,1]$ 
if for any $\bar{a}$, $\bar{b}$ from ${\bf X}^n$, 
$$ 
|F(\bar{a} )-F(\bar{b})| \le \delta (d(\bar{a} ,\bar{b} )). 
$$ 
\end{definicja}  

The choice of the connectives above guarantees that 
the following statement holds (see \cite{IMI-APAL}). 

\begin{lem} \label{ContMod} 
For any continuous relational structure $M$, 
where each $n$-ary relation has $n\cdot \mathsf{id}$ as an inverse 
continuity modulus, any formula $\phi$ admits an inverse continuity 
modulus which is of the form $k \cdot \mathsf{id}$, where $k$ 
depends on the complexity of $\phi$. 
\end{lem}

\begin{remark} \label{rel} 
{\em By Lemma 4.1 of \cite{BNT} each $n$-ary functional symbol $F$ can be replaced 
by the predicate  $D_F (\bar{x},y) = d(F(\bar{x}), y)$. 
It is clear that the continuity moduli with respect to variables from $\bar{x}$ are the same 
and $\mathsf{id}$ works as a continuity modulus for $y$. 
Thus we may always assume that $L$ is relational. }
\end{remark} 

For a continuous structure $M$ defined on $({\bf Y},d)$ let 
$\mathsf{Aut}(M)$ be the subgroup of $\mathsf{Iso} ({\bf Y})$ consisting of 
all isometries preserving the values of atomic formulas. 
It is easy to see that $\mathsf{Aut}(M)$ is a closed subgroup with 
respect to the topology on $\mathsf{Iso} ({\bf Y})$ defined above. 

For every $c_1 ,...,c_n \in M$ and $A\subseteq M$ 
we define the $n$-type $\mathsf{tp}(\bar{c}/A)$ of $\bar{c}$ over $A$ 
as the set of all $\bar{x}$-conditions with parameters from $A$ 
which are satisfied by $\bar{c}$ in $M$.  
Let $S_n (T_A )$ be the set of all $n$-types over $A$ 
of the expansion of the theory $T$ by constants from $A$. 
There are two natural topologies on this set. 
The {\bf logic topology} is defined by the basis consisting of 
sets of types of the form $[\phi (\bar{x})<\varepsilon ]$, 
i.e. types containing some $\phi (\bar{x})\le \varepsilon'$ with 
$\varepsilon '<\varepsilon$.   
The logic topology is compact. 

The $d$-topology is defined by the metric 
$$
d(p,q)= \mathsf{inf} \{  d(\bar{c} ,\bar{b})| \mbox{ there is a model } M \mbox{ with } M\models p(\bar{c})\wedge q(\bar{b})\}. 
$$ 
By Propositions 8.7 and 8.8 of \cite{BYBHU} the $d$-topology is finer 
than the logic topology and $(S_n (T_A ),d)$ is a complete space.

The following notion is helpful when we study some concrete 
examples, for example the Urysohn space. 
A relational continuous structure $M$ is 
{\bf approximately ultrahomogeneous} if for any $n$-tuples 
$(a_1 ,..,a_n )$ and $(b_1 ,...,b_n )$ with the same 
quantifier-free type (i.e. with the same values of 
predicates for corresponding subtuples) and any 
$\varepsilon >0$ there exists $g\in \mathsf{Aut}(M)$ such that 
$$ 
\mathsf{max} \{ d(g(a_j ),b_j ): 1\le j \le n\} \le \varepsilon . 
$$  
As we already mentioned any Polish group can be chosen 
as the automorphism group of a continuous metric 
structure which is approximately ultrahomogeneous.  

The bounded Urysohn space $\mathfrak{U}$ (see Section 2.3) 
is {\bf ultrahomogeneous} in the traditional sense: 
any partial isomorphism between two tuples extends 
to an automorphism of the structure \cite{U}.  
Note that this obviously implies that $\mathfrak{U}$ is 
approximately ultrahomogeneous. 

We will use the continuous version of 
$L_{\omega_1 \omega}$ from \cite{BYIov} (see also \cite{BNT}. 
We remind the reader that continuous 
$L_{\omega_1 \omega}$-formulas are defined by the standard 
procedure applied to countable conjunctions and disjunctions 
(see \cite{BYIov}). 
Each continuous infinite formula depends on finitely many free variables. 
The main demand is the existence of continuity moduli of such 
formulas.  
It is usually assumed that a continuity modulus 
$\delta_{\phi ,x}$ satisfies the equality 
$$ 
\delta_{\phi ,x} (\varepsilon) = \mathsf{sup} \{ \delta_{\phi ,x} (\varepsilon '): 0<\varepsilon' <\varepsilon \} 
$$ 
and 
$$ 
\delta_{\bigwedge \Phi ,x}(\varepsilon )= \mathsf{sup} \{ \delta'_{\bigwedge \Phi ,x}(\varepsilon '): 0<\varepsilon' <\varepsilon \},  
\mbox{ where } \delta'_{\bigwedge \Phi ,x}= \mathsf{inf} \{ \delta_{\phi ,x}: \phi\in \Phi \} .  
$$ 

\subsection{Logic action} \label{LA}

Fix a countable continuous signature 
$$
L=\{ d,R_1 ,...,R_k ,..., F_1 ,..., F_l ,...\}
$$ 
and a Polish space $({\bf Y},d)$. 
Let $S$ be a dense countable subset of ${\bf Y}$. 
Let $\mathsf{seq}(S)=\{ \bar{s}_i :i\in \omega\}$ be the set (and an enumeration) 
of all finite sequences (tuples) from $S$. 
Let us define the space of metric $L$-structures on $({\bf Y},d)$. 
Using the recipe as in the case of $\mathsf{Iso}({\bf Y})$ 
we introduce a metric on the set of $L$-structures as follows. 
Enumerate all tuples of the form $(\varepsilon ,j,\bar{s})$, where  
$\varepsilon \in \{ 0,1\}$ and when $\varepsilon =0$,  $\bar{s}$ is a tuple 
from $\mathsf{seq}(S)$ of the length of the arity of $R_j$,  and for $\varepsilon =1$, 
$\bar{s}$ is a tuple from $\mathsf{seq}(S)$ of the length of the arity of $F_j$.  
For metric $L$-structures $M$ and $N$ let   
$$ 
\delta_{\mathsf{seq}(S)} (M ,N )= \sum_{i=1}^{\infty} \{ 2^{-i} |R^M_j (\bar{s} )-R^N_j (\bar{s} )| 
\footnote{resp. $ 2^{-i} d(F^M_j (\bar{s} ),F^N_j (\bar{s} ))$ when $\varepsilon =1$}  
\mbox{ : }  i \mbox{ is the number of } (\varepsilon ,j,\bar{s}) \} . 
$$ 
Since the predicates and functions are uniformly continuous 
(with respect to moduli of $L$) and $S$ is dense 
in ${\bf Y}$, we see that $\delta_{\mathsf{seq}(S)}$ is a complete metric. 
Moreover by an appropriate choice of rational values for 
$R_j (\bar{s})$ we find a countable dense subset of metric 
structures on ${\bf Y}$, i.e. the space obtained is Polish. 
We denote it by ${\bf Y}_L$. 
It is clear that $\mathsf{Iso} ({\bf Y})$ acts on ${\bf Y}_L$ continuously. 
Thus we consider ${\bf Y}_L$ as an $\mathsf{Iso}({\bf Y})$-space and 
call it the {\bf space of the logic action} on ${\bf Y}$. 

\begin{remark} \label{present}
{\em It is worth noting that in this definition a structure on ${\bf Y}$ 
(say $M$) is identified with its presentation $\mathsf{Diag} (M,S)$, 
see Introduction.}  
\end{remark}

It is convenient to consider the following {\bf basis} 
of the topology of ${\bf Y}_L$. 
Fix a finite sublanguage $L'\subset L$, a finite subset 
$S'\subset S$, a finite tuple $q_1 ,...,q_t \in {\bf Q} \cap [0,1]$ 
and a rational $\varepsilon \in [0,1]$ with $1-\varepsilon < 1/2$.  
Consider a diagram $D$ of $L'$ on $S'$ of some inequalities of the form   
$$
d(F_j (\bar{s}) ,s' ) > \varepsilon \mbox{ , } d(F_j (\bar{s}) ,s' )< 1 - \varepsilon , 
$$
$$ 
|R_j (\bar{s}) - q_i | > \varepsilon \mbox{ , } |R_j (\bar{s}) - q_i |< 1 - \varepsilon ,\mbox{ with } \bar{s}\in \mathsf{seq}(S'), s'\in S' . 
$$
(i.e. in the case of relations we consider 
negations of statements of the form: 
$|R_j (\bar{s}) - q_i |\le \varepsilon$ , 
$|R_j (\bar{s}) - q_i |\ge 1 - \varepsilon$).  
The set of metric $L$-structures realizing $D$ is 
an open set of the  topology of ${\bf Y}_L$ and 
the family of sets of this form is a basis of this topology.   
Compactness theorem for continuous logic (see \cite{BYU}) 
shows that the topology is compact. 
We will call it {\bf logic} too. 

If in Remark \ref{present} one relax the conditions on formulas 
$\phi$ used in $\mathsf{Diag} (M,S)$ (for example allowing 
$\phi$ to be from some $L_{\omega_1 \omega}$-fragments) 
the topology can become richer (and the basis should be corrected).  
Moreover by the continuous version of the Lopez-Escobar theorem 
(\cite{BNT}, \cite{CL}) every Polish group action arises 
as an action of some closed $G\le \mathsf{Iso}({\bf Y})$  
on the space of separable continuous structures of some 
$L_{\omega_1\omega}$-sentence, \cite{CL}. 
This possibility will be discussed in the next section. 

\section{Good and nice topologies} 

In Section 2.2 we give the main concepts of 
the generalized model theory. 
They are based on the notion of a grey subset
introduced in \cite{BYM}. 
The corresponding preliminaries are given in Section 2.1. 
In Section 2.3 we describe the most important examples
of the situation. 
In Section 2.4 we demonstrate several 
applications of our approach. 
They concern complexity of some subsets of the logic space. 
In a sense this section explains the reason why the questions of complexity 
are considered under the framework of good/nice topologies (of Section 2.2). 

\subsection{Grey subsets} 

The notion of grey subsets was introduced in \cite{BYM}. 
It has become very fruitful, see \cite{BNT}, \cite{CL}, \cite{IMI-APAL}
and \cite{Chen}. 

A function $\phi$ from a space ${\bf X}$ to $[-\infty ,+\infty ]$ 
is {\bf upper (lower) semi-continuous} if the set 
$\phi_{<r}$ (resp. $\phi_{>r}$) is open for all $r\in \mathbb{R}$ 
(here $\phi_{<r} = \{ z\in {\bf X}: \phi (z) <r \}$, a {\bf cone}). 
A {\bf grey subset} of ${\bf X}$, denoted 
$\phi \sqsubseteq {\bf X}$, is a function 
${\bf X}\rightarrow [0,\infty ]$. 
It is {\bf open (closed)}, $\phi \sqsubseteq_o {\bf X}$ 
(resp. $\phi \sqsubseteq_c {\bf X}$), 
if it is upper (lower) semi-continuous. 
We also write $\phi \in {\bf \Sigma}_1$ when 
$\phi \sqsubseteq_o {\bf X}$ and we write 
$\phi \in {\bf \Pi}_1$ when $\phi \sqsubseteq_c {\bf X}$.  
We will assume below that values of 
a grey subset belong to $[0,1]$. 

\bigskip 

Let us return to the situation of Section \ref{LA}. 
We fix a language $L$, a countable dense subset $S$ of ${\bf Y}$ 
and study subsets of  ${\bf Y}_L$. 
One of the basic observations is that any first-order continuous sentence  
$\phi (\bar{c})$, $\bar{c}\in S$, defines a grey subset of ${\bf Y}_L$: 
$$
\phi (\bar{c}) \mbox{ takes } M \mbox{ to the value } \phi^{M} (\bar{c}). 
$$ 
Moreover Proposition \ref{EsLo} below says that   
$\phi (\bar{c})$ defines a grey subset of ${\bf Y}_L$ which belongs 
to ${\bf \Sigma}_n$  for some $n$.  
It is Proposition 1.1 in \cite{IMI-APAL}. 

\begin{prop} \label{EsLo}
For any continuous formula $\phi(\bar{v})$ of the 
language $L$ there is a natural number $n$ such that 
for any tuple $\bar{a}\in S$ 
and $\varepsilon \in [0,1]$, the subset 
$$
Mod(\phi ,\bar{a},<\varepsilon )=\{ M :M \models\phi (\bar{a})<\varepsilon \} 
$$
$$
\mbox{ ( or } 
Mod(\phi ,\bar{a},>\varepsilon )=\{ M :M \models\phi (\bar{a})>\varepsilon \} \mbox{ ) } 
$$ 
of the space ${\bf Y}_L$ of $L$-structures, belongs to ${\bf \Sigma}_n$. 
\end{prop}

When $G$ is a Polish group, then a grey subset 
$H \sqsubseteq G$ is called a {\bf grey subgroup} if 
$$
H(1)=0 \mbox{ , }\forall g\in G (H(g)=H(g^{-1})) \mbox{ and } 
\forall g, g'\in G (H(gg')\le H(g)+H(g')). 
$$ 
This is equivalent to Definition 2.5 from \cite{BYM}. 
It is worth noting that by Lemma 2.6 of \cite{BYM} 
an open grey subgroup is clopen. 

\bigskip

If $H$ is a grey subgroup, then for every $g\in G$ we define 
the grey coset $Hg$ and the grey conjugate $H^g$ as follows:
$$
\begin{array}{l@{\ = \ }l}
Hg(h)&H(hg^{-1})\\
H^g(h)&H(ghg^{-1}).
\end{array}
$$
Observe that if $H$ is open, then $Hg$ is an open grey subset and
$H^g$ is an open grey subgroup.

\begin{definicja} 
Let ${\bf X}$ be a continuous $G$-space.  
A grey subset $\phi \sqsubseteq {\bf X}$ is called {\bf invariant} with 
respect to a grey subgroup $H \sqsubseteq G$ if 
for any $g\in G$ and $x\in {\bf X}$ we have $\phi (g(x)) \le \phi (x) \dot+ H(g)$. 
\end{definicja}  

Since $H(g)=H(g^{-1})$, the inequality from the definition is 
equivalent to $\phi (x) \le \phi (g(x)) \dot+ H(g)$. 

\bigskip 

\begin{remark} {\em (see Section 2.1 of \cite{IMI-APAL}).
It is clear that for every continuous structure $M$ 
(defined on ${\bf Y}$) any continuous formula $\phi (\bar{x})$ 
defines a clopen grey subset of $M^{|\bar{x}|}$. 
Moreover note that when $\phi (\bar{x},\bar{c})$ is a continuous 
formula with parameters $\bar{c}\in M$ and $\delta$ is a linear inverse 
continuous modulus for $\phi (\bar{x}, \bar{y})$ (see Definition \ref{inverse}), 
then $\phi$ is invariant with respect to the open grey subgroup 
$H_{\delta, \bar{c}} \sqsubseteq \mathsf{Aut}(M)$ defined by  
$$
H_{\delta, \bar{c}}(g) =\delta ( d((c_1 ,\ldots ,c_n ), (g(c_1 )),\ldots ,g(c_n ))) \mbox{, where } g\in \mathsf{Aut}(M) ,
$$
i.e.  
$$
\phi(g(\bar{a}),\bar{c}) \le \phi (\bar{a},\bar{c}) + H_{\delta ,\bar{c}}(g) .    
$$ 
}  
\end{remark} 
\bigskip 
 
In the space of continuous $L$-structures $ {\bf Y}_L$ this 
remark has the follows version (see Lemma 2.2 in \cite{IMI-APAL}). 

\begin{lem} 
Let $\delta$ be an inverse continuity modulus for $\phi (\bar{x})$, 
which is linear. 
The grey subset defined by $\phi (\bar{c}) \sqsubseteq {\bf Y}_L$ 
is invariant with respect to the {\bf grey stabiliser}  
$H_{\delta ,\bar{c}}\sqsubseteq \mathsf{Iso}({\bf Y})$ defined 
as follows. 
$$
H_{\delta ,\bar{c}}(g) =\delta ( d((c_1 ,\ldots ,c_n ), (g(c_1 )),\ldots ,g(c_n ))) \mbox{, where } g\in \mathsf{Iso}({\bf Y}) . 
$$
\end{lem}

\subsection{Nice bases}

In this section we consider a certain class of Polish $G$-spaces. 
To describe it we need the following definition. 
 
\begin{definicja}   \label{gbasis} 
A family $\mathcal{U}$  of open grey subsets of a Polish space 
${\bf X}$ with a topology $\tau$ is called a {\bf grey basis}  
of $\tau$ if the family 
$\{\phi_{<r}:\phi\in{\mathcal{U}}, r\in {\mathbb{Q}}\cap (0,1)\}$ 
is a basis of $\tau$.
\end{definicja} 

We now describe our typical assumptions on $G$:  
\begin{itemize} 
\item $G$ is a Polish group; 
\item we distinguish a countable dense subgroup $G_0 <G$ and 
a countable family of clopen grey subsets $\mathcal{R}$ of $G$ 
which is a grey basis of the topology of $G$; 
\item we assume that $\mathcal{R}$ consists of all  
$G_0$-cosets of grey subgroups from $\mathcal{R}$, i.e. 
for each $\rho\in\mathcal{R}$ there is a grey subgroup 
$H\in \mathcal{R}$ and an element $g_0 \in G_0$ 
so that for any $g\in G$, $\rho (g ) = H(g g^{-1}_0 )$; 
\item  we assume that $\mathcal{R}$ is closed under 
$G_0$-conjugacy, under $\mathsf{max}$ and truncated 
multiplication by positive rational numbers. 
\end{itemize}

\begin{remark}\label{Gbasis}  
{\em In Remark 2.9 of \cite{IMI-APAL}
it is observed that for every Polish group $G$ there is a a countable $G_0 <G$ and 
a countable family of open grey subsets $\mathcal{R}$ 
satisfying these assumptions. }
\end{remark} 

We will see below that if the space $({\bf Y},d)$ is good enough 
(for example the bounded Urysohn space) and 
$S$ is a dense countable subset of ${\bf Y}$, 
then  the family $\mathcal{R}$ of grey subsets of $G=\mathsf{Iso}({\bf Y})$ 
can be chosen among  grey cosets of the form 
$$ 
\rho (g) = q \cdot d(\bar{b}, g(\bar{a})) \mbox{ , where } q\in \mathbb{Q}^{+} 
\mbox{ and } 
$$ 
$$ 
\bar{a}, \bar{b} \mbox{ are tuples from } S 
\mbox{ which are isometric in } {\bf Y}.  
$$ 
If the metric is bounded by 1 we mean the truncated multiplication by $q$ 
in the formula above. 

When we fix $G_0$, $\mathcal{R}$ and consider 
a Polish $G$-space $({\bf X},d)$
we also distinguish a countable grey basis  
$\mathcal{U}$ of the topology of ${\bf X}$.  
Let $\tau$ be 
the corresponding topology.  

The approach of  generalized model theory of H. Becker from \cite{becker2} 
suggests that along with the $d$-topology $\tau$ we shall consider
some special topology on ${\bf X}$ which is called {\bf nice}.
In the case of Polish $G$-spaces this idea has been 
realized in \cite{IMI-APAL} with using continuous logic.   
Since we do not need the corresponding material in 
exact form we introduce the following very general definition. 

\begin{definicja}  \label{NB}
Let $\mathcal{R}$ be a grey basis of $G$ consisting of cosets of 
open grey subgroups of $G$ which also belong to $\mathcal{R}$. 
Assume that the subfamily of $\mathcal{R}$ of all open grey subgroups 
is closed under $\mathsf{max}$ and truncated multiplication by numbers from 
$\mathbb{Q}^{+}$. 

We say that a family $\mathcal{B}$ of Borel grey subsets
of the $G$-space $(\mathbf{X}, \tau )$ is a {\bf good basis} 
with respect to $\mathcal{R}$ if: \\ 
(i) $\mathcal{B}$ is countable and generates the topology finer than $\tau$;\\ 
(ii) for each $\phi\in \mathcal{B}$ there exists an open grey 
subgroup $H\in \mathcal{R}$ such that $\phi$ is $H$-invariant.
\end{definicja}

It will be usually assumed that 
all constant functions $q$, $q\in \mathbb{Q}\cap [0,1]$ 
are in $\mathcal{B}$.

\begin{definicja} \label{nto}
A topology ${\bf t}$ on ${\bf X}$ is $\mathcal{R}$-{\bf good} for the $G$-space
$\langle {\bf X}, \tau\rangle$ if the following conditions are satisfied.\\
(a) The topology ${\bf t}$ is Polish, ${\bf t}$ is finer than $\tau$
and the $G$-action remains continuous with respect to ${\bf t}$. \\
(b) There exists a grey basis $\mathcal{B}$ of ${\bf t}$ which is good with respect to $\mathcal{R}$. 
\end{definicja}

Nice bases and nice topologies introduced in \cite{IMI-APAL} are good. 
We remind the reader that a good basis $\mathcal{B}$ with respect to $\mathcal{R}$
is {\bf nice} if the following additional properties hold: 
\begin{quote}  
(iii) for all $\phi_1, \phi_2 \in \mathcal{B}$, the functions $\neg \phi_1$, 
$\mathsf{min}(\phi_1 ,\phi_2)$, $\mathsf{max}(\phi_1 ,\phi_2)$, $|\phi_1 - \phi_2 |$, 
$\phi_1 \dot- \phi_2$ $\phi_1 \dot+ \phi_2$ belong to $\mathcal{B}$;\\ 
(iv) for all $\phi \in \mathcal{B}$ and $q\in \mathbb{Q}^{+}$   
the truncated product $q\cdot \phi$   
belongs to $\mathcal{B}$; \\ 
(v) for all $\phi\in \mathcal{B}$ and open grey subsets  
$\rho \in \mathcal{R}$ the Vaught transforms (see Section 2.1 in \cite{IMI-APAL}) 
$\phi^{*\rho}, \phi^{\Delta \rho}$ belong $\mathcal{B}$.  
\end{quote} 

In the situation of standard examples (see Section 2.3) 
the property that the basis is good is straightforward. 
It is much more difficult to prove that the basis is nice. 
Theorem 3.2 from \cite{IMI-APAL} is an example of a result of this kind. 
The following theorem gives existence of nice topologies. 
This is Theorem 2.12 in \cite{IMI-APAL}. 
Note that the assumptions 
on the grey basis $\mathcal{R}$ 
follow from the conditions of/before  Remark \ref{Gbasis}.

\begin{thm} \label{existence} 
Let $G$ be a Polish group and  $\mathcal{R}$ be a countable grey basis satisfying 
the assumptions of Definition \ref{NB} and  the following closure property: 
\begin{quote} 
for every grey subgroup $H\in\mathcal{R}$ and every $g\in G$ 
if $Hg\in\mathcal{R}$, then $H^g\in\mathcal{R}$. 
\end{quote} 
Let $\langle {\bf  X}, \tau\rangle $ be a Polish $G$-space and $\mathcal{F}$ 
be a countable family of Borel grey subsets of ${\bf X}$ generating 
a topology finer than $\tau$ such that for any $\phi\in \mathcal{F}$  
there is a grey subgroup $H\in \mathcal{R}$ such that $\phi$ 
is invariant with respect to $H$. 
 
Then there is an $\mathcal{R}$-nice topology for $G$-space 
$\langle {\bf  X}, \tau\rangle$ 
such that $\mathcal{F}$ consists of open grey subsets.
\end{thm}

In Section 2.4 we describe possible applications of statements of this kind.

\subsection{Countable approximating substructures} \label{CAS}

In this section we give basic examples of good bases 
and topologies on some logic spaces.  
The following definition is taken from \cite{BYBM}. 
\begin{definicja} 
Let $({\bf M},d)$ be a Polish metric structure with universe ${\bf M}$. 
We say that a (classical) countable structure
$N$ is a {\bf countable approximating substructure} 
of ${\bf M}$ if the following conditions are
satisfied: 
\begin{itemize}  
\item The universe $N$ of $N$ is a dense countable subset of
$({\bf M}, d)$.  
\item Any automorphism of $N$ extends to a (necessarily unique) 
automorphism of ${\bf M}$, and $\mathsf{Aut}(N)$ is dense in
$\mathsf{Aut}({\bf M})$.
\end{itemize} 
\end{definicja} 

Let $G_0$ be a dense countable subgroup of $\mathsf{Aut}(N)$.  
We may consider it as a subgroup of $\mathsf{Aut} ({\bf M})$. 

{\bf Family $\mathcal{R}^{{\bf M}}(G_0 )$.} 
Let $\mathcal{R}_0$ be the family of all clopen 
grey subgroups of $\mathsf{Aut}({\bf M})$ of the (truncated) form 
$$ 
H_{q, \bar{s}} : g\rightarrow q \cdot d(g(\bar{s}), \bar{s}), 
\mbox{ where } \bar{s}\subset N, \mbox{ and } q\in \mathbb{Q}^{+}.   
$$ 
It is clear that $\mathcal{R}_0$ is closed under 
conjugacy by elements of $G_0$. 
Consider the closure of $\mathcal{R}_0$ under 
the function $\mathsf{max}$ and define $\mathcal{R}^{{\bf M}}(G_0 )$ 
to be the family of all $G_0$-cosets of grey 
subgroups from $\mathsf{max}(\mathcal{R}_0 )$. 
Then $\mathcal{R}^{{\bf M}} (G_0 )$ is countable and the 
family of all $(H_{q, \bar{s}})_{<l}$ where 
$H\in \mathcal{R}_0$ and $l\in \mathbb{Q}$, 
generates the topology of $\mathsf{Aut}({\bf M},d)$.  
Moreover it is easy to see that $G_0$ and 
$\mathcal{R}^{{\bf M}} (G_0 )$ satisfy all the conditions of/before  
Remark \ref{Gbasis} for $\mathcal{R}$ and in particular $\mathcal{R}^{{\bf M}} (G_0 )$ 
satisfies the conditions of Theorem \ref{existence} for $\mathcal{R}$. 

{\bf Family $\mathcal{B}_{\mathcal{L}}$. } 
Let $L$ be a relational language of a continuous signature with 
inverse continuity moduli $\le n\cdot \mathsf{id}$ for $n$-ary relations. 
We will assume that $L$ extends the language of the structure ${\bf M}$. 

Let $\mathcal{L}$ be a countable fragment of $L_{\omega_1 \omega}$, 
in particular $\mathcal{L}$ be closed under first-order connectives. 
Note that inverse continuity moduli of first-order continuous formulas 
(with connectives as in Introduction) can be taken linear 
(of the form $k\cdot \mathsf{id}(x)$). 
Thus it is easy to see that  every formula of $\mathcal{L}$ 
has linear inverse continuity moduli. 

Let $\mathcal{B}_{\mathcal{L}}$ be the family of all grey subsets defined 
on the logic space ${\bf M}_L$ by continuous $\mathcal{L}$-sentences (with parameters) 
as follows  
$$
\phi (\bar{s}) : M\rightarrow \phi^{M} (\bar{s}), 
\mbox{ where }\bar{s}\in N \mbox{ and } \phi (\bar{x}) \in \mathcal{L} . 
$$ 
By linearity of inverse continuity moduli it is easy to see 
that for any continuous sentence $\phi (\bar{s})$ 
there is a number $q\in \mathbb{Q}$ 
(depending on the continuity modulus of $\phi$) such that 
the grey subset as above is $H_{q, \bar{s}}$-invariant. 
As a result we have the following statement. 

\begin{quote} 
Let $\mathcal{B}_{\mathcal{L}}$ be a family of grey subsets corresponding 
to a countable continuous fragment $\mathcal{L}$ of $L_{\omega_1 \omega}$. 
Then the family $\mathcal{B}_{\mathcal{L}}$ is a 
good basis with respect to $\mathcal{R}^{{\bf M}} (G_0 )$. 
\end{quote} 
\bigskip 

{\bf (A)} 
Let us consider the following example. 
The Urysohn space of diameter 1  
is the unique Polish metric space of diameter 1 
which is universal and ultrahomogeneous. 
This space $\mathfrak{U}$ is considered in 
the continuous signature  $\langle d \rangle$. 

The countable counterpart of $\mathfrak{U}$ is the 
{\bf rational Urysohn space of diameter 1}, $\mathbb{Q}\mathfrak{U}$, 
which is both ultrahomogeneous and universal for countable 
metric spaces with rational distances and diameter $\le 1$. 
The space $\mathfrak{U}$ is interpreted as ${\bf M}$ above and 
$\mathbb{Q}\mathfrak{U}$ will be our $N$. 
It is shown in Section 6.1 of \cite{BYBM} that there is 
an embedding of $\mathbb{Q}\mathfrak{U}$ into $\mathfrak{U}$ so that: \\
(i) $\mathbb{Q}\mathfrak{U}$ is an approximating substructure of $\mathfrak{U}$: 
it is dense in $\mathfrak{U}$; 
any isometry of  $\mathbb{Q}\mathfrak{U}$ extends to an isometry of  
$\mathfrak{U}$ and $\mathsf{Iso} (\mathbb{Q}\mathfrak{U})$ is dense in $\mathsf{Iso}(\mathfrak{U})$;  \\ 
(ii) for any $\varepsilon>0$, any partial isometry $h$ of  
$\mathbb{Q}\mathfrak{U}$ with domain $\{ a_1 ,...,a_n\}$ and any isometry 
$g$ of $\mathfrak{U}$ such that $d(g(a_i ),h(a_i ))<\varepsilon$ 
for all $i$, there is an isometry $\hat{h}$ of  $\mathbb{Q}\mathfrak{U}$ 
that extends $h$ and is such that for all  $x\in \mathfrak{U}$, 
$d(\hat{h}(x),g(x))<\varepsilon$.  

Let $G_0$ be a dense countable subgroup of $\mathsf{Iso}(\mathbb{Q}\mathfrak{U})$.  
By (i) we may consider it as a subgroup of $\mathsf{Iso} (\mathfrak{U})$. 
We now define $\mathcal{R}^{\mathfrak{U}} (G_0 )$ by the recipe above.  
As we already know $G_0$ and 
$\mathcal{R}^{\mathfrak{U}} (G_0 )$ satisfy all the conditions of/before  
Remark \ref{Gbasis} and in particular $\mathcal{R}^{\mathfrak{U}} (G_0 )$ 
satisfies the conditions of Theorem \ref{existence}. 

\bigskip 
 
Let $L$ be a relational language of a continuous signature as above. 
Let $\mathcal{L}$ be a countable fragment of $L_{\omega_1 \omega}$ and  
let $\mathcal{B}_{\mathcal{L}}$ be the family of all grey subsets defined 
by continuous $\mathcal{L}$-sentences (with parameters from $\mathbb{Q}\mathfrak{U}$) as above. 
We already know that $\mathcal{B}_{\mathcal{L}}$ is a good basis. 
It is proved in \cite{IMI-APAL} (see Theorem 3.2) that this basis is nice. 

\begin{thm} \label{mainUrysohn} 
The family $\mathcal{B}_{\mathcal{L}}$ is a 
$\mathcal{R}^{\mathfrak{U}}(G_0 )$-nice basis. 
\end{thm}

\bigskip

{\bf (B) A separable Hilbert space.} 
We follow \cite{BYBM} and \cite{ros09}. 
Let us consider the complex Hilbert space $l_2(\mathbb{N})$. 
Let ${\cal Q}$ denote the algebraic closure of
$\mathbb{Q}$, and consider the countable subset
${\cal Q}l_2$ of $l_2 (\mathbb{N})$ of all sequences with finite 
support and coordinates from ${\cal Q}$. 
It is shown in Section 6.2 of \cite{BYBM} (with use Section 7 of \cite{ros09}), 
that it is an approximating substructure of $l_2 (\mathbb{N})$. 
In particular we have another pair playing the role of $({\bf M},N)$. 
Since $l_2 (\mathbb{N})$ is unbounded, the authors of \cite{BYBM} consider
instead its closed unit ball, equipped with functions
$x \rightarrow \alpha x$ for $|\alpha | \le 1$ 
and
$(x, y)\rightarrow \frac{x+y}{2}$, from which
$l_2 (\mathbb{N})$  can be recovered. 
According Remark \ref{rel} we will consider a relational language for 
this structure. 
 
The automorphism group of the unit ball is ${\bf U}(l_2 (\mathbb{N}))$, 
the unitary group of the whole complex Hilbert space $l_2(\mathbb{N})$.  
The topology of pointwise convergence is the strong 
operator topology.  

Let $G_0$ be a dense countable subgroup of ${\bf U}(\mathcal{Q}l_2 )$.  
We may consider it as a subgroup of ${\bf U}(l_2 (\mathbb{N}))$. 
We now apply the procedure of 
$\mathcal{R}^M (G_0 )$ and $\mathcal{B}_{\mathcal{L}}$. 
As a result we obtain the 
{\bf family $\mathcal{R}^{{\bf H}}(G_0 )$} and 
a grey basis defined on ${\bf U}(l_2 (\mathbb{N}))$. 

Let $L$ be a relational language of a continuous signature extending the language 
of the unit ball and satisfying the assumptions above and  
let $\mathcal{L}$ be a countable fragment of $L_{\omega_1 \omega}$. 
Let $\mathcal{B}_{\mathcal{L}}$ be the corresponding family of all grey subsets 
of the logic space $l_2 (\mathbb{N})_{L}$. 
This is a good basis with respect to 
$\mathcal{R}^{{\bf H}}(G_0 )$. 

\bigskip 

{\bf (C) The measure algebra on $[0,1]$.}
Denote by $\lambda$ 
the Lebesgue measure on the unit interval
$[0,1]$. 
We view its automorphism group
$\mathsf{Aut}([0,1], \lambda )$ as the automorphism group of the
Polish metric structure
$$
( MALG,0,1,\wedge , \vee , \neg , d), 
$$ 
where MALG denotes the measure algebra on $[0,1]$ 
and
$d(A, B) = \lambda (A \Delta B)$
(see [Kec95]). 
The approximating substructure is the countable
measure algebra $A$ generated by dyadic intervals. 
This is observed in Section 6.3 od \cite{BYBM}. 
Exactly as in the case of $\mathfrak{U}$ and $l_2 (\mathbb{N})$ 
one can define a family of open grey subgroups 
of $\mathsf{Aut}([0,1], \lambda )$, say $\mathcal{R}^{Aut}(G_0)$,  
and a good bases of the corresponding logic spaces 
with respect to $\mathcal{R}^{Aut}(G_0 )$.

\begin{remark} 
{\em The basic continuous metric structures which appear in {\bf (A)} - {\bf (C)}, 
i.e. $\mathfrak{U}$, the unit ball of $l_2 (\mathbb{N})$ and MALG, 
are ultahomogeneous structures in the classical sense: any partial isomorphism 
between two tuples extends to an automorphism of the structure. 
This is in particular mentioned in Section 3.1 of \cite{BYFra}. }
\end{remark}

\subsection{The Effros Borel structure of $\mathsf{Iso} (\mathfrak{U})$. Applications} 

Given a Polish space ${\bf Y}$ let $\mathcal{F} ({\bf Y})$ 
denote the set of closed subsets of ${\bf Y}$. 
The Effros structure on $\mathcal{F} ({\bf Y})$ 
is the Borel space with respect to the $\sigma$-algebra generated by 
the sets 
$$ 
\mathcal{C}_U = \{ D\in \mathcal{F} ({\bf Y}): D\cap U \not=\emptyset \}, 
$$ 
for open $U \subseteq {\bf Y}$. 
For various ${\bf Y}$ this space serves for analysis 
of Borel complexity of families of closed subsets 
(see \cite{KNT} and \cite{RZ} some recent results). 
It is convenient to use the fact that there is a sequence 
of Kuratowski-Ryll-Nardzewski selectors 
$s_n : \mathcal{F} ({\bf Y}) \rightarrow {\bf Y}$, $n\in \omega$, 
which are Borel functions such that for every non-empty 
$F\in \mathcal{F} ({\bf Y})$ the set $\{ s_n (F) ; n\in \omega \}$ 
is dense in $F$. 

Given a Polish group $G$ and a continuous (or Borel) 
action of $G$ on a Polish space ${\bf Y}$ one can consider 
the Borel space 
$$
\mathcal{F} ({\bf Y} )^m \times \mathcal{F} (G)^n. 
$$
In the situation when ${\bf Y}$ and $G$ have  grey bases
$\mathcal{B}$ and $\mathcal{R}$ respectively 
which satisfy the conditions of Theorem \ref{existence}   
one can consider ${\bf Y}^m$ with respect to the good topology induced by $\mathcal{B}$ (say ${\bf t}$). 
Then 
many natural Borel subsets of ${\bf Y}^m \times G^n$ can be viewed as elements of  
$$
\mathcal{F} (({\bf Y}, {\bf t}) )^m \times \mathcal{F} (G)^n. 
$$

By Theorem 2.2 of \cite{CL} for any Polish group $G$ 
and any standard Borel $G$-space ${\bf X}$
there is a continuous group monomorphism 
$\Phi : G \rightarrow \mathsf{Iso} (\mathfrak{U} )$ 
and a Borel $\Phi$-equivariant injection $f: {\bf X} \rightarrow \mathfrak{U}_L$. 
We only need here that the language $L$ is countable relational 
with 1-Lipschitz symbols of unbounded arity.  
As a result all Polish groups can be considered as elements of 
$\mathcal{F} (\mathsf{Iso} (\mathfrak{U} ))$, all Polish spaces 
are elements of $\mathcal{F} (\mathfrak{U}_L)$ and Polish $G$-spaces 
are pairs from 
$$
\mathcal{F} (\mathfrak{U}_L ) \times \mathcal{F} (\mathsf{Iso} (\mathfrak{U} )). 
$$ 
Let $\mathcal{R}^{\mathfrak{U}}(G_0 )$ and $\mathcal{B}_{\mathcal{L}}$ be grey bases defined in Section \ref{CAS} in the case of 
$\mathsf{Iso}( \mathfrak{U} )$ and $\mathfrak{U}$.  
Let ${\bf t}$ be the corresponding nice topology. 
The following proposition is a version of a well-known fact. 

\begin{prop} \label{elem-Borel} 
(1) The following relations  from  $(\mathcal{F} (\mathsf{Iso} (\mathfrak{U})))^2$,  
$(\mathcal{F} (\mathfrak{U}_L ))^2$,  $(\mathcal{F} (\mathfrak{U}_L ,{\bf t}))^2$,
$(\mathcal{F} (\mathsf{Iso} (\mathfrak{U} )))^3$, 
$\mathcal{F} (\mathsf{Iso} (\mathfrak{U} ))\times \mathcal{F} (\mathfrak{U}_L ) \times \mathcal{F} (\mathfrak{U}_L )$ and 
$\mathcal{F} (\mathsf{Iso} (\mathfrak{U} ))\times \mathcal{F} (\mathfrak{U}_L ,{\bf t}) \times \mathcal{F} (\mathfrak{U}_L, {\bf t} )$ 
(under natural interpretations) are Borel: 
$$
\{ (A,B): A \subseteq B \} \mbox{ , }  \{ (A,B,C): AB \subseteq C\}. 
$$ 
(2) The closed subgroups of $\mathsf{Iso} (\mathfrak{U})$ 
form a Borel set $\mathcal{U} (\mathsf{Iso} (\mathfrak{U}))$ in 
$\mathcal{F} (\mathsf{Iso} (\mathfrak{U}))$. \\
(3) The Polish $G$-spaces form a Borel set in 
$ \mathcal{F} (\mathsf{Iso} (\mathfrak{U} ))\times \mathcal{F} (\mathfrak{U}_L ) \times\mathcal{F} (\mathfrak{U}_L )$ 
and closed $G$-subspaces of $(\mathfrak{U}_L ,{\bf t})$ form a Borel set in 
$ \mathcal{F} (\mathsf{Iso} (\mathfrak{U} ))\times\mathcal{F} (\mathfrak{U}_L ,{\bf t})\times\mathcal{F} (\mathfrak{U}_L ,{\bf t})$ . 
\end{prop} 

{\em Proof.} Statement (2) is well-known: see Section 3.2 of \cite{RZ}.  
Moreover statements (1) and (2) are variants of Lemmas 2.4 and 2.5 from \cite{KNT} 
which were proved for  $S_{\infty}$. 
It is also  mentioned in \cite{KNT} that they hold in general. 
Statement (3) follows  from (1) and (2). 
We only mention here that a Polish $G$-space is viewed as a triple consisting of $G$, 
the subspace and the graph of the action. 
$\Box$ 

\bigskip 

In model theory a theory $T_1$ is a model companion of $T$ 
if $T_1$ is model complete and every model of $T$ embeds into a model of $T_1$ 
and vice versa. 
One of the definitions of model completeness  states that any formula 
is equivalent to an existential one (or a universal one). 

In the case of the logic space $\mathfrak{U}_L$ theories are identified with 
${\bf t}$-closed invariant subsets. 
It is convenient to fix an enumeration of the sets 
$$ 
\mathcal{B}_{\mathcal{L}}(\mathbb{Q}) = \{ (\phi )_{< r} : 
\phi \in \mathcal{B}_{\mathcal{L}} \mbox{ and } r\in \mathbb{Q} \cap [0,1] \} ,    
$$ 
$$
\mathcal{B}_{o\mathcal{L}}(\mathbb{Q}) = \{ (\phi )_{< r} : \phi 
\mbox{ is a clopen member of } \mathcal{B}_{\mathcal{L}} \mbox{ and } 
r\in \mathbb{Q} \cap [0,1] \} ,    
$$ 
where the latter one is a basis of the topology $\tau$. 
The following definition is a version of Definition 1.3 from \cite{MI05} 
and Proposition 1.4 of \cite{IMI-APAL}. 
 
\begin{definicja} 
Let $X_0$ and $X_1$ be closed invariant subsets of  $(\mathfrak{U}_L ,{\bf t})$. 
We say that $X_1$ is a {\bf companion} of $X_0$ if $\tau$-closures of 
$X_0$ and $X_1$ coincide and any element of $\mathcal{B}_{\mathcal{L}}$ 
is $\tau$-clopen on $X_1$. 
\end{definicja} 

\begin{thm} \label{comp}
The set of pairs $(X_0, X_1 )$ of 
$\mathsf{Iso} (\mathfrak{U} )$-invariant members of 
$\mathcal{F} (\mathfrak{U}_L ,{\bf t})$ 
with the condition that $X_1$ is a companion of $X_0$ 
is Borel. 
\end{thm} 

{\em Proof.} 
Applying Proposition \ref{elem-Borel} we consider pairs $(X_0 ,X_1 )$  
of ${\bf t}$-closed invariant subsets 
as elements of the corresponding Borel set of triples 
$(\mathsf{Iso} (\mathfrak{U} ), X_0 , X_1 )$.  
Using Kuratowski-Ryll-Nardzewski selectors the condition that $\tau$-closures 
of $X_0$ and $X_1$ are the same can be written as follows:  
$$ 
(\forall A_k \in \mathcal{B}_{o\mathcal{L}}(\mathbb{Q}))\forall i \exists j \exists l ( 
s_i (X_1 ) \in A_k \rightarrow s_j (X_0 ) \in A_k ) \wedge 
( s_i (X_0 ) \in A_k \rightarrow s_l (X_1 ) \in A_k ). 
$$ 
Now note that for any two ${\bf t}$-closed $A$ and $B$ 
the condition 
$A\cap X_1 \subseteq B\cap X_1$ 
is equivalent to the formula: 
$$ 
\forall j (s_j (X_1 )\in A \rightarrow s_j (X_1 ) \in B ). 
$$ 
In particular this condition is Borel. 
We can now express that any element of $\mathcal{B}_{\mathcal{L}}(\mathbb{Q})$ 
is $\tau$-clopen on $X_1$ as follows: 
$$ 
\forall i (\forall B_l \in \mathcal{B}_{\mathcal{L}}(\mathbb{Q})) 
(\exists A_k \in \mathcal{B}_{o\mathcal{L}}(\mathbb{Q})) 
(\exists A_m \in \mathcal{B}_{o\mathcal{L}}(\mathbb{Q})) 
(s_i (X_1 ) \in B_l \rightarrow s_i (X_1) \in A_k \wedge 
$$
$$
(A_k \cap X_1\subseteq B_l \cap X_1)) \wedge 
(s_i (X_1 ) \not\in B_l \rightarrow s_i (X_1) \in A_m )\wedge 
(A_m \cap X_1 \subseteq X_1 \setminus B_l )). 
$$
$\Box$ 

\bigskip 

\begin{remark} 
{\em The theorem above is a counterpart of the statement that identifying theories 
a language $L$ with closed subsets of the compact space of complete $L$-theories 
the binary relation to be a model companion is Borel. 
Although the authors have not found it in literature, it is true and 
possibly is folklore. }
\end{remark} 
 
Theorem \ref{comp} confirms that the approach of good/nice topologies is useful.  
It provides  a topological tool for a general property from logic (model companions).

\section{Computable presentations}

If ${\bf X}$ is of the form ${\bf Y}_L$ then it makes sense 
to study complexity of sets of indices of computable structures 
of natural model-theoretic classes. 
In the case of first order structures this approach is 
traditional, see \cite{AK}, \cite{EG}, \cite{EG1} and \cite{GK}. 
In Section 3.1 we give an appropriate generalization and 
in Section 3.2 we illustrate it in the case of 
$\mathfrak{U}_L$ for relational $L$. 

These ideas were already presented by the authors in Section 5 
of preprint \cite{IMI-old}. 
We have discovered that they are  closely related to the approach of 
the recent paper of A.G. Melnikov and A. Montalb\'{a}n \cite{MelMon}. 
In fact the main concern below is to realize the situation of Sections 2.1 - 2.2 of 
\cite{MelMon} in cases {\bf (A)} - {\bf (C)}  of Section 2.3. 
Having this we arrive in a field where the results of \cite{MelMon} work. 
 
It is worth noting here that the approach of Section 3.1 can be applied when 
one considers computable presentations of Polish spaces, 
see \cite{mos} and \cite{W}. 
Some details are given in Remark \ref{mosc}. 
In a sense this is the easiest case. 
The approach of \cite{MelMon} also works here.

\subsection{Computable grey subsets}

Consider the situation of Section 2.2. 
Let $G$ be a Polish group and $\mathcal{R}$ 
be a distinguished countable family of clopen 
grey cosets which is a grey basis of $G$: 
\begin{itemize}  
\item the family 
$\{ \rho_{<q} : \rho\in \mathcal{R}$ and $q\in \mathbb{Q}^{+}\cap [0,1]\}$ 
forms a basis of the topology of $G$.   
\end{itemize} 
We fix a countable dense subgroup $G_0 <G$ so that:
\begin{itemize} 
\item  $\mathcal{R}$ is closed under $G_0$-conjugacy and consists of all  
$G_0$-cosets of grey subgroups from $\mathcal{R}$;   
\item the set of grey subgroups from 
$\mathcal{R}$ is closed under $\mathsf{max}$ and truncated multiplication 
by positive rational numbers.
\end{itemize} 
Let  
$$ 
\mathcal{R}^+_{\mathbb{Q}} =\{ \sigma_{<r} :\sigma \in \mathcal{R}\mbox{ and }r\in \mathbb{Q}^+ \cap [0,1]\} , 
$$ 
$$ 
\mathcal{R}^-_{\mathbb{Q}}=\{ \rho_{>q} : \rho\in \mathcal{R}\mbox{ and }q\in \mathbb{Q}^{+}\cap [0,1]\} . 
$$ 
We assume that 
\begin{itemize} 
\item 
there is a computable 1-1-enumeration of the family 
$\mathcal{R}_{\mathbb{Q}} = \mathcal{R}^+_{\mathbb{Q}} \cup \mathcal{R}^-_{\mathbb{Q}}$ 
so that the relation of inclusion between members 
of this family is computable. 
\end{itemize} 

\begin{remark} 
{\em Having this assumption we arive at the case 
that $(G,\mathcal{R}_{\mathbb{Q}})$ is a 
{\bf computably presented $\omega$-continuous domain},  
see \cite{EH} and \cite{ES}. 
In fact our assumptions are slightly stronger. 
Moreover in  {\bf A1} - {\bf A4} below 
we will make them much stronger.}
\end{remark} 

\begin{itemize} 
\item Let $({\bf X},\tau )$ be  a Polish $(G,G_0 ,\mathcal{R})$-space 
together with a distinguished countable $G_0$-invariant 
grey basis $\mathcal{U}$ (see Definition \ref{gbasis}) 
of clopen grey subsets which is closed under $\mathsf{max}$ 
and truncated multiplying by positive rational numbers.   
\item Let 
$$ 
\mathcal{U}^+_{\mathbb{Q}} =\{ \sigma_{<r} :\sigma \in \mathcal{U}\mbox{ and }r\in \mathbb{Q}^+ \cap [0,1]\} \mbox{ (a basis of } ({\bf X}, \tau ) \mbox{ ), }  
$$ 
$$ 
\mathcal{U}^-_{\mathbb{Q}}=\{ \rho_{>q} : \rho\in \mathcal{U}\mbox{ and }q\in \mathbb{Q}^{+}\cap [0,1]\} \mbox{ and } 
\mathcal{U}_{\mathbb{Q}} = \mathcal{U}^+_{\mathbb{Q}} \cup \mathcal{U}^-_{\mathbb{Q}} ,  
$$ 
and the relation of inclusion between sets of $\mathcal{U}_{\mathbb{Q}}$
be computable (under an appropriate computable coding). 
\end{itemize} 
As a result $({\bf X}, \mathcal{U}_{\mathbb{Q}})$
is a computably presented $\omega$-continuous domain.

Note that in the discrete case these circumstances are 
standard and in particular arise when one studies 
computability in $S_{\infty}$-spaces of logic actions.

\begin{remark} \label{mosc}
{\em It is worth noting that when we have a recursively 
presented Polish space in the sense of the book of 
Moschovakis \cite{mos} (Section 3), 
then a basis of the form $\mathcal{U}_{\mathbb{Q}}$ 
as above (in fact $\mathcal{U}^+_{\mathbb{Q}}$) 
can be naturally defined. 
Indeed, let us recall that a recursive presentation of a Polish space 
$({\bf X},d)$ is any sequence $S_{{\bf X}} =\{ x_i : i\in\omega\}$ 
which is a dense subset of ${\bf X}$ satisfying 
the condition that $(i,j,m,k)$-relations 
$$
d(x_i ,x_j )\le \frac{m}{k+1} \mbox{ and } d(x_i ,x_j ) < \frac{m}{k+1} 
$$ 
are recursive. 
If in this case for all $i$ we define grey subsets 
$\phi_i (x) = d(x,x_i )$,  then all balls 
$(\phi_i )_{<r}$, $r\in \mathbb{Q}$, form a basis 
$\mathcal{U}^+_{\mathbb{Q}}= \{ B_i :i\in\omega \}$ of ${\bf X}$ 
which under appropriate enumeration 
(together with co-balls $(\phi_i )_{>r}$) 
satisfies our requirements above. 
When $G$ is a Polish group with a left-invariant 
metric $d$, then for any $q_1 ,...,q_k\in \mathbb{Q}$ 
and any tuple $h_1 ,...,h_k \in G$ the grey subset 
$\phi_{\bar{q},\bar{h}} (x)= \mathsf{max}_{i\le k} (q_i \cdot d(h_i , xh_i ))$ 
is a grey 
subgroup\footnote{apply $d(h_i ,xyh_i ) \le d(h_i ,x h_i ) + d(xh_i ,xy h_i )= d(h_i ,xh_i )+d(h_i ,yh_i )$ 
together with the fact that $max$ applied to grey subgroups gives grey subgroups again}. 
If $G$ is a recursively presented space 
with respect to a dense countable subgroup 
$G_0$ and the multiplication is recursive, 
then let $\mathcal{V}$ consist of all $\phi_{\bar{q},\bar{h}}$ 
with $\bar{h} \in G_0$ and let $\mathcal{R}$ consist 
of all $G_0$-cosets of these grey subgroups. 
The structure (domain) $(\mathcal{R}_{\mathbb{Q}}, \subset )$ 
is computably presented. 

If $G$ isometrically acts on ${\bf X}$ and $x_1 ,...,x_k$ 
is a finite subset of the recursive presentation $S_{{\bf X}}$ 
then as we already know the function 
$\psi _{\bar{q},\bar{x}}(g)=\mathsf{max}_{i\le k} (q_i \cdot d(x_i ,g(x_i )))$ 
also defines a grey subgroup. 
When $G$ has a recursive multiplication and a recursive action 
on ${\bf X}$ (see Section 3 of \cite{mos}) so that 
the recursive presentation $S_{{\bf X}}$ is $G_0$-invariant, 
then let $\mathcal{V}$ consist of these subgroups 
and $\mathcal{R}$ consist of the $G_0$-cosets.  
Then the structure $(\mathcal{R}_{\mathbb{Q}}, \subset )$ 
is computably presented. } 
\end{remark} 

As we will see below the following assumptions are satisfied in 
the majority of interesting cases.  

\paragraph{ Computability assumptions.}  

{\bf A1.} We assume that under our 1-1-enumerations of the families 
$\mathcal{R}_{\mathbb{Q}}$ and $\mathcal{U}_{\mathbb{Q}}$ 
the sets of indices of $\mathcal{U}^+_{\mathbb{Q}}$, $\mathcal{R}^+_{\mathbb{Q}}$  
and the set of rational cones 
$$
\mathcal{V}^+_{\mathbb{Q}} = \{ H_{<r} : H \mbox{ is a graded subgroup from } \mathcal{R} \} 
$$ 
$$
\mathcal{V}^-_{\mathbb{Q}} = \{ H_{>r} : H \mbox{ is a graded subgroup from } \mathcal{R} \}
$$
are distinguished by computable unary relations on $\omega$. 

{\bf A2.} We assume that under our 1-1-enumerations of the families 
$\mathcal{R}_{\mathbb{Q}}$ and $\mathcal{U}_{\mathbb{Q}}$ 
the binary relation to be in the pair $\sigma_{<r}$, $\sigma_{>r}$ 
for $\sigma \in \mathcal{R}$ or $\sigma\in \mathcal{U}$ 
is computable. 

{\bf A3.} We also assume that the following relation is computable: 
\begin{quote} 
{\it 
$Inv(V,U)\Leftrightarrow (V\in \mathcal{V}^+_{\mathbb{Q}})\wedge (U\in \mathcal{U}^+_{\mathbb{Q}})\wedge (U$ 
is $V$-invariant $)$ }.
\end{quote} 
By invariantness we mean the property that $U$ is presented 
as $\{ x \in {\bf X}: \phi (x)<r \}$, the set $V$ is presented as $\{ g\in G : H(g)<s\}$ 
and $\phi$ is an $H$-invariant grey subset 
(in particular the inequality $\phi (g(x)) < r+s$ holds for $x\in U$ and $g\in V$). 

{\bf A4.}  We assume that there is an algorithm deciding 
the problem whether for a natural number $i$ and for a basic set 
of the form $\sigma_{<r}$ for $\sigma$ 
from $\mathcal{U}$ or $\mathcal{R}$ and $r\in \mathbb{Q}$,  
the diameter of $\sigma_{<r}$  is less than $2^{-i}$.  

\bigskip

Under this setting we introduce the main 
notion, which is a counterpart of  a {\bf computable structure}. 
In \cite{mos} in the case of a recursively presented 
Polish space it is defined that a point $x\in {\bf X}$ 
is {\bf recursive} if the set 
$\{ s :x\in B_s , B_s \in \mathcal{U}_{\mathbb{Q}}\}$ 
is computable. 
We imitate it in the following definition.

\begin{definicja} \label{Comp} 
We say that an element $x\in {\bf X}$ is {\bf computable} if the relation 
$$
Sat_x (U)\Leftrightarrow (U \in \mathcal{U}_{\mathbb{Q}})\wedge (x\in U)
$$ 
is computable. 
\end{definicja}

In the case of the logic action of $S_{\infty}$, 
when $x$ is a structure on $\omega$ and all $H$ and $\phi$ are 
two-valued, this notion is obviously equivalent to 
the notion of a computable structure.  

We will denote by $Sat_x ({\mathcal{U}}_{\mathbb{Q}})$ the set 
$\{ C: C\in {\mathcal{U}}_{\mathbb{Q}}$ and $Sat_x (C)$ holds $\}$. 

\begin{remark} 
{\em In \cite{MelMon} computable topological spaces are 
considered under so called } formal inclusion $\ll$    
{\em (it corresponds to terms "approximation" or "way-below" in other 
sources). 
In Definition 2.2 of \cite{MelMon} it is defined for 
computable Polish metric spaces, 
but in fact this relation can be defined in more general situations. 
Axioms (F1) - (F4) given in \cite{MelMon} after Definition 2.2 
describe the field of applications of this notion. 
It is always assumed in \cite{MelMon} that $\ll$ is computably enumerable. 
In our framework this relation can be defined as follows: 
$$ 
\sigma_{<r } \ll \sigma'_{<r'} \Leftrightarrow \exists r_1 (r<r_1 \wedge \sigma_{<r_1} \subseteq \sigma'_{<r'} 
\wedge \mathsf{diam} (\sigma_{<r}) \le \frac{1}{2} \mathsf{diam}(\sigma'_{<r'})). 
$$ 
Then it is computably enumerable and satisfies (F1) - (F4) of \cite{MelMon}. 
In particular the results of \cite{MelMon} hold in the cases of Section 3.2 below.    
We only add here that in \cite{MelMon} computable elements are those $x$ 
for which  $Sat_x ({\mathcal{U}}_{\mathbb{Q}})$ is computably enumerable. }
\end{remark}

We now make few basic observations which are very helpful when one tries 
to estimate the complexity of some families of computable structures.  
The following lemma follows  from the assumption 
that $\mathcal{U}$ is a grey basis and satisfies {\bf A4}.   

\begin{lem} \label{kappa} 
If $x\in {\bf X}$ is computable then there is a computable function 
$\kappa :\omega \rightarrow \mathcal{U}^+_{\mathbb{Q}}$ such that for all 
natural numbers $n$, $x\in \kappa (n)$ and $diam(\kappa (n))\le 2^{-n}$. 
\end{lem} 
 
We also say that 
\begin{quote} 
an element $g\in G$ is {\bf computable} 
if the relation $(N\in \mathcal{R}_{\mathbb{Q}})\wedge (g\in N)$ 
is computable. 
\end{quote} 
Then there is a computable function realizing the 
same property as $\kappa$ above but already in the case 
of the basis $\mathcal{R}_{\mathbb{Q}}$. 

In the following lemma we use standard indexations 
of the set of computable functions and of the set 
of all finite subsets of $\omega$. 

\begin{lem} \label{21}
The following relations belong to $\Pi^0_2$:\\
(1) $\{ e:$ the function $\varphi_e$ is a characteristic 
function of a subset of $\mathcal{U}_{\mathbb{Q}} \}$; \\ 
(2) $\{ (e,e'):$ there is a computable element $x\in {\bf X}$ 
such that the function $\varphi_e$ is a characteristic function 
of the set $Sat_x ({\mathcal{U}_{\mathbb{Q}}})$ and 
the function $\varphi_{e'}$ realizes the corresponding function 
$\kappa$ defined in Lemma \ref{kappa} $\}$;\\   
(3) $\{ (e,e'):$ there is an element $g\in G$ such that the 
function $\varphi_e$ is a characteristic function of the 
subset $\{ N\in\mathcal{R}_{\mathbb{Q}}:g\in N\}$ and the function 
$\varphi_{e'}$ realizes the corresponding function $\kappa$ 
defined as in Lemma \ref{kappa} $\}$.    
\end{lem} 

{\em Proof.} (1) Obvious. 
Here and below we use the fact that a function is computable if 
and only if its graph is computably enumerable. 

(2) Under {\bf A1} and {\bf A4} the corresponding definition can be described as follows: 
$$
("e \mbox{ is a characteristic function of a subset of }\mathcal{U}_{\mathbb{Q}}")\wedge 
$$ 
$$
(\forall n)((\varphi_{e'}(n)\in\mathcal{U}^+_{\mathbb{Q}})\wedge (\varphi_{e'}(n)\not=\emptyset )
\wedge (\varphi_{e}(\varphi_{e'}(n))=1 )\wedge 
(diam(\varphi_{e'}(n))<2^{-n} )) \wedge  
$$
\begin{quote} 
$(\forall d)(\exists n )(($ 
"every element $U'$ of the finite subset of $\mathcal{U}_{\mathbb{Q}}$ 
with the canonical index $d$ satisfies $\varphi_e (U') =1$") 
$\leftrightarrow ($ "$\varphi_{e'}(n)$ is contained in any element $U'$ 
of the finite subset of $\mathcal{U}_{\mathbb{Q}}$ with the canonical index $d$"$))$.  
\end{quote}
The last part of the conjunction ensures that the intersection of 
any finite subfamily of $\mathcal{U}_{\mathbb{Q}}$ of cones $U'$ 
with $\varphi_e (U') =1$ contains a closed cone of the form $\phi_{\le r}$ 
of sufficiently small diameter. 
Now the existence of the corresponding $x$ follows by 
Cantor's intersection theorem for complete spaces. 

(3) is similar to (2). 
$\Box$
\bigskip 

We say that $e$ is an index of a computable 
element $x\in {\bf X}$ if $\varphi_e$ is a characteristic function of 
$Sat_x (\mathcal{U}_{{\bf Q}})$. 
We now have the following straightforward proposition. 

\begin{prop} The set of indices of computable 
elements of ${\bf X}$ belongs to $\Sigma^0_3$. 
\end{prop} 

\subsection{Computable approximating structures}

In this section we show that the computability assumptions 
{\bf A1} - {\bf A4} given in Section 3.1 are satisfied in the case of 
good graded bases presented in {\bf (A) } - {\bf (C)} 
of Section 2.3. 
In fact we only consider case {\bf (A)}. 
Cases {\bf (B)},{\bf (C)} are similar.

\paragraph{(D) Computable presentation of the logic space over $\mathfrak{U}$. }
Let $L$ be a relational language satisfying assumptions of Section 2.3. 
Let us consider the space $\mathfrak{U}_L$ and the family of 
grey cosets $\mathcal{R}^{\mathfrak{U}}(G_0 )$ defined in Section 2.3 {\bf (A)}. 
The latter will be interpreted as $\mathcal{R}$ of Section 3.1. 

To define the {\bf  grey basis} $\mathcal{U}$  of Section 3.1 we use 
the recipe of the definition of the basis of the topology of ${\bf Y}_L$ in Section 1.3. 
For a finite sublanguage $L'\subset L$, a finite subset 
$S'\subset \mathbb{Q}\mathfrak{U}$ and a finite tuple 
$q_1 ,...,q_t \in \mathbb{Q} \cap [0,1]$  
consider the maximum $\mathsf{max}$ of some grey subsets of the form 
$$ 
|R_j (\bar{s}) - q_i |  \mbox{ , } 1 \dot{-} |R_j (\bar{s}) - q_i | ,\mbox{ with } \bar{s}\in \mathsf{seq}(S'), s'\in S' . 
$$  
When $\sigma$ is this maximum, the  inequality $\sigma <\varepsilon$ 
corresponds to a basic open set of the  topology of ${\bf Y}_L$ 
as in Section 1.3.

Let $\mathcal{L}$ be the fragment of all first order continuous formulas. 
Let $\mathcal{B}_0$ be the nice basis corresponding to $\mathcal{L}$ 
(see Theorem \ref{mainUrysohn}).  
It is worth noting that the {\bf grey basis} 
$\mathcal{U}$ is a subfamily 
of the family of all grey subsets from $\mathcal{B}_0$.  
Moreover $\mathcal{U}$ corresponds to quantifier free 
$L$-formulas.  

To verify that the $\mathsf{Iso} (\mathfrak{U})$-space 
$\mathfrak{U}_L$ and the bases $\mathcal{R}$, $\mathcal{U}$ satisfy 
the computability conditions of Section 3.1 (in particular {\bf A1} - {\bf A4}), 
we need the following proposition. 

\begin{prop} \label{QU-decidability} 
The elementary theory of the structure 
$\mathbb{Q}\mathfrak{U}$ in the binary language of inequalities 
$$ 
d(x,x') \le (\mbox{ or  } \ge ) q \mbox{ , where } q\in \mathbb{Q}\cap [0,1], 
$$ 
extended by all constants from $\mathbb{Q}\mathfrak{U}$ is decidable. 
\end{prop} 

{\em Proof.} 
It is well known that $\mathbb{Q}$ can be identified with the natural numbers so that 
the ordering of the rational numbers becomes a computable relation. 
Thus the language in the formulation can be considered as a computable one. 

It is noticed in \cite{KPT} that the first order structure $\mathbb{Q}\mathfrak{U}$ is 
universal ultrahomogeneous in the language 
$$
d(x,x')= q \mbox{ , where } q\in \mathbb{Q}\cap [0,1].  
$$ 
This obviously implies that $\mathbb{Q}\mathfrak{U}$ is a universal 
ultrahomogeneous first-order structure in the language of inequalities 
as in the statement of the proposition.  
So one can present $\mathbb{Q}\mathfrak{U}$ in this language as a Fra\"{i}ss\'{e} limit 
of an effective sequence of finite structures.  
Enumerating elements of structures from this sequence and 
describing distances between them, we obtain an effective set 
of axioms of the form 
$$ 
d(c,c') \le (\mbox{ or } \ge ) q \mbox{ , where } c,c' \in \mathbb{Q}\mathfrak{U} \mbox{ and } q\in \mathbb{Q}\cap [0,1].  
$$ 
We also add all standard $\forall \exists$-axioms 
stating that the age of $\mathbb{Q}\mathfrak{U}$ is an amalgamation class. 
The obtained axiomatization describes a complete theory having elimination of quantifiers. 
$\Box$ 

\bigskip 

\begin{cor} \label{presentation}
The structure $\mathbb{Q}\mathfrak{U}$ under the language of binary relations 
$$ 
d(x,y) \le (\mbox{ or } \ge ) q \mbox{ , where } q\in \mathbb{Q}\cap [0,1], 
$$ 
has a presentation on $\omega$ 
so that all relations first-order definable 
in $\mathbb{Q}\mathfrak{U}$, are decidable. 
\end{cor} 

This obviously follows from 
Proposition \ref{QU-decidability}. 
Let us fix such a presentation. 

\paragraph{Coding $\mathcal{R}_{\mathbb{Q}}$, cones of grey cosets.} 
Let 
$$ 
H_{q, \bar{s}} : g\rightarrow q \cdot d(g(\bar{s}), \bar{s}), 
\mbox{ where } \bar{s}\subset \mathbb{Q}\mathfrak{U}, \mbox{ and } q\in \mathbb{Q}^{+}.   
$$ 
be a grey subgroup and $g_0 \in G_0$ take $\bar{s}'$ to $\bar{s}$.  

Then we can code the $*q'$-cone of the grey coset  
$$
H_{q ,\bar{s}} g_0 :  
g\rightarrow q \cdot d(g(\bar{s}'), \bar{s}), 
$$ 
by the number of the tuple $(q, \bar{s},\bar{s}',q',*)$, 
where $\bar{s}$, $\bar{s}'$ are identified with the corresponding 
tuples from $\omega$ with respect to the presentation of 
Corollary \ref{presentation} and 
$*$ is one of the symbols $<, \le ,>, \ge$.  
Note that the tuples $\bar{s}$, $\bar{s}'$ have the same 
quantifier free diagram (which is determined by a finite subdiagram). 
By Corollary \ref{presentation} the set of all 
tuple $(q, \bar{s},\bar{s}',q',*)$ of this form is computable 
and by ultrahomogeneity of the structure from this corollary 
they code all possible cones. 

To see that  the relation of inclusion between cones 
of this form is decidable note that 
\begin{quote} 
$(q, \bar{s},\bar{s}',q',*)$ defines a subset of 
the cone of $(q_1 , \bar{s}_1 ,\bar{s}'_1 ,q'_1 ,*_1 )$ if 
for every tuple $\bar{s}''\bar{s}''_1$ of the same quantifier free type with $\bar{s}'\bar{s}'_1$ 
which also satisfies the $*$-inequality between $q\cdot d(\bar{s}'',\bar{s})$ and $q'$, 
the corresponding $*_1$-inequality   between $q_1\cdot d(\bar{s}''_1,\bar{s}_1)$ and $q'_1$ 
holds.  
\end{quote} 
Indeed by ultrahomogeneity this exactly states that if for an automorphism 
$g$  the $*$-inequality between $q\cdot d(g(\bar{s}'),\bar{s})$ and $q'$ holds, 
then the corresponding $*_1$-inequality   between $q_1\cdot d(g(\bar{s}'_1),\bar{s}_1)$ and $q'_1$ 
also holds. 
Thus to decide the inclusion problem between these cones it suffices to 
formulate the statement above as a formula (with parameters  $\bar{s}',\bar{s}$, $\bar{s}'_1,\bar{s}_1)$) 
and to verify if it holds in the structure $\mathbb{Q}\mathfrak{U}$.

{\bf Cones of grey subgroups}  
(i.e. the set $\mathcal{V}_{\mathbb{Q}}$)
are distinguished in the set of codes of $\mathcal{R}_{\mathbb{Q}}$ 
by the computable subset of tuples as above 
with $\bar{s}=\bar{s'}$. 
\parskip0pt 

\paragraph{Coding $\mathcal{U}_{\mathbb{Q}}$.} 
Since we interpret elements of $\mathcal{B}_{0}$ by 
first order $L$-formulas with parameters from $\mathbb{Q}\mathfrak{U}$ 
and without free variables,
it is obvious that both $\mathcal{B}_{0}$ and $\mathcal{U}$ 
can be coded in $\omega$ so that the operations of 
connectives are defined by computable functions. 
Moreover $\mathcal{U}$ is a decidable subset of $\mathcal{B}_0$. 
Thus the elements of the grey basis $\mathcal{U}$ are coded as a computable set. 
Now all cones of the form 
$\sigma_{<q}$, $\sigma_{>q}$, $\sigma_{\le q}$, $\sigma_{\ge q}$ 
can be enumerated so that all natural relations between them 
(in particular relations from {\bf A2}) are computable. 
For example if $S'$ is a finite subset of $\mathbb{Q}\mathfrak{U}$ and cones 
$\sigma_{ <q}$ and $\sigma'_{<q'}$ correspond 
to inequalities of the form 
$$ 
|R_j (\bar{s}) - q_i |   <q  \mbox{  , } 1 \dot{-} |R_j (\bar{s}) - q_i | <q   
  \mbox{  ,  with } \bar{s}\in \mathsf{seq}(S') 
 \mbox{  (and similarly in the case of } \sigma'_{<q' }\mbox{  ), }   
$$  
then it can be verified effectively if the inequalities of the cone $\sigma_{< q}$ 
follow from the ones of the cone $\sigma'_{<q'}$ together with the diagram of the metric 
on $S'$ and the inequalities provided by the continuity moduli.   
It is worth noting here that the diagram of the metric on $S'$ 
(i.e. all equalities $d(s,s') = q$ with $s,s'\in S'$, $q\in [0,1]\cap \mathbb{Q}$) 
is decidable by Corollary \ref{presentation}. 

\paragraph{Satisfying A3.}  
Let $U$ be of the form $\sigma_{<q}$ for 
$\sigma\in \mathcal{U}$ and $V$ 
be of the form $H_{<k}$ for $H\in \mathcal{V}$.  
We assume that the inequalities of $\sigma_{<q}$ are as in the previous paragraph. 
Let $S'$ be a finite subset of $\mathbb{Q}\mathfrak{U}$ which contains all parameters 
which appear in the definition of $U$ and $H$.  
Since $\mathbb{Q}\mathfrak{U}$ is ultrahomogeneous the condition $Inv (V,U)$ is satisfied if 
and only if the following property holds.  
\begin{quote} 
If $\gamma$ is a partial isometry of $\mathbb{Q}\mathfrak{U}$ with domain $S'$ 
and fixing the parameters appearing in $H$, then  
the value of $\sigma$ (with respect to the parameters $S'$) 
is preserved under $\gamma$.  
\end{quote} 
When this property holds for all possible interpretations of 
the $L$-symbols on $\mathfrak{U}$ it just follows from the diagram of 
the metric on $S'$ and the inequalities provided by the continuity moduli. 
Thus by Corollary \ref{presentation} this relation is  decidable. 
 
\paragraph{Satisfying A4.} 
Let $\sigma$ be a $\mathsf{max}$-formula of the previous paragraphs  
which defines an element of $\mathcal{U}$.  
To compute $\mathsf{diam}(\sigma_{<q})$ consider the definition of 
the metric $\delta_{\mathsf{seq}(\mathbb{Q}\mathfrak{U})}$ of the space 
$\mathfrak{U}_{L}$ with respect to $\mathsf{sec}(\mathbb{Q}\mathfrak{U})$ 
in the beginning of Section 1.1. 
Find all numbers $i$ of 
tuples $(j,\bar{s}')$ such that 
$R_j (\bar{s}')$  appears in $\sigma$.   
We may assume that appearance of such subformulas  
forces inequalities of the form 
$q'_i \le R_j (\bar{s}') \le q_i$ for rational 
$0\le q'_i <q_i \le 1$. 
Let $I$ be the (finite) subset of such $i$. 
Then $\mathsf{diam}(\sigma_{<q})$ is computed by  
$$ 
\sum_{i=1}^{\infty} \{ 2^{-i}\mbox{ : }  i \not\in I\} + 
\sum_{i\in I} 2^{-i} |q_i - q'_i| . 
$$  
In particular we have an algorithm for comparing it with powers $2^{-i}$. 

The case of basic clopen sets of $\mathcal{R}^U$ is similar. 

\bigskip

\paragraph{ (E) Decidability.} 

We have found that our arguments for the computability assumptions 
{\bf A1} - {\bf A4} can be applied 
for some other related questions and possibly the most natural 
ones are involved into decidability of continuous theories.   
This explains why we now consider this issue. 

Let $\Gamma$ be a set of continuous formulas 
of a continuous signature $L$ with a metric. 
Let $\phi$ be a continuous $L$-formula. 

\begin{definicja} 
(see \cite{BYP}, Section 9) 
The value    
$\mathsf{sup} \{ \phi^{M} :  \mbox{ for } M \models \Gamma =0\}$
is called the {\bf degree of truth  of} $\phi$ {\em with respect to} $\Gamma$.  
We denote this value by $\phi^{\circ}$.  
\end{definicja} 

If the language $L$ is computable, the set of 
all continuous $L$-formulas and the set of all $L$-conditions of the form 
$$ 
\phi \le \frac{m}{n} \mbox{ , where } \frac{m}{n}\in \mathbb{Q}_+ , 
$$ 
are computable. 

We remind the reader that a real number $r\ge 0$ is {\bf computable} 
if there is an algorithm which for any natural number $n$ 
finds a natural number $k$ such that 
$$ 
 \frac{k-1}{n} \le r \le \frac{k+1}{n}.  
$$ 

Corollary 9.11 of \cite{BYP} states that when $\Gamma$ 
is computably enumerable 
and $\Gamma=0$ axiomatizes a complete theory, 
then the value $\phi^{\circ}$ is 
a recursive real which is  uniformly computable from $\phi$. 
The latter exactly means that the corresponding 
complete theory is {\bf decidable}. 
Note that in this case the value $\phi^{\circ}$  
coincides with the value of $\phi$ in models of $\Gamma =0$.   

The following theorem shows that in the situations 
of examples of Section 2.3 the expansion of the structure 
by the countable approximating substructure  
has decidable continuous theory.

\begin{thm} 
The structure $(\mathfrak{U},  s)_{ s\in\mathbb{Q}\mathfrak{U}}$ of 
the expansion of the bounded Urysohn space  by constants from 
$\mathbb{Q}\mathfrak{U}$ has decidable continuous theory. 

The same statement holds for structures $({\bf M}, s)_{s\in N}$ 
where $({\bf M},d) \in \{ l_2 (\mathbb{N}), MALG \}$ and $N$ 
is the corresponding countable approximating substructure
see Section 2.3, {\bf (B)} and {\bf (C)}).  
\end{thm} 

{\em Proof.} 
To prove the theorem we use Corollary 9.11 of \cite{BYP}. 
We only consider the case of $(\mathfrak{U},  s)_{ s\in\mathbb{Q}\mathfrak{U}}$. 
The remaining cases are similar. 
 
Let $T_{\mathbb{Q}\mathfrak{U}}$ be the set of the standard axioms of $\mathfrak{U}$ 
(with rational $\varepsilon$ and $\delta$, see Section 5 in \cite{U}) 
together with all quantifier free axioms 
describing distances between constants from $\mathbb{Q}\mathfrak{U}$. 
We claim that the set $T_{\mathbb{Q}\mathfrak{U}}$ is computable. 
Since the set of all standard axioms of $\mathfrak{U}$ 
is computable (see \cite{U}), it suffices to check that 
the set of all axioms of the form 
$$ 
d(c,c') = q \mbox{ , where } c,c' \in \mathbb{Q}\mathfrak{U} \mbox{ and } q\in \mathbb{Q}\cap [0,1], 
$$ 
is computable. 
This follows from the fact that 
the elementary (not continuous) theory of the structure 
$\mathbb{Q}\mathfrak{U}$ in the language of binary relations 
together with all constants $c\in \mathbb{Q}\mathfrak{U}$ is decidable, 
Proposition \ref{QU-decidability}. 

Note that $T_{\mathbb{Q}\mathfrak{U}}$ axiomatizes 
the continuous theory of a single continuous structure, 
i.e. the corresponding continuous theory is complete.  
Indeed, otherwise there is a separable continuous structure 
$M\models T_{\mathbb{Q}\mathfrak{U}}$ such that for some  
tuple $\bar{s}\in  \mathbb{Q}\mathfrak{U}$ the structures $(\mathfrak{U},\bar{s})$ 
and the reduct of $M$, say $M'$, to the signature $(d,\bar{s})$, 
do not satisfy the same inequalities of  the form 
$$
\phi (\bar{s})\le (<) q \mbox{ or }  \phi (\bar{s})\ge (>) q \mbox{ where } 
q \in \mathbb{Q}\cap [0,1].  
$$
On the other hand since $\mathfrak{U}$ is separably  categorical 
(see Section 4) and ultrahomogeneous, the structures $M'$ and  $(\mathfrak{U},\bar{s})$ 
are isomorphic, contradicting the previous sentence.  

By Corollaries 9.8 and  9.11 of \cite{BYP} there is an algorithm which for every 
continuous sentence $\phi(\bar{s})$ computes its value in $\mathfrak{U}$. 
$\Box$

\bigskip

\begin{remark} 
{\em It is worth noting that when we apply Proposition \ref{QU-decidability} 
we only need computability of the set of axioms of the form  
$$ 
d(c,c') = q \mbox{ , where } c,c' \in \mathbb{Q}\mathfrak{U} \mbox{ and } q\in \mathbb{Q}\cap [0,1], 
$$ 
This can be shown as in the proof of Proposition \ref{QU-decidability}. 
Moreover the corresponding argument works 
in the cases of $l_2 (\mathbb{N})$ and $MALG$. } 
\end{remark}

\section{Complexity of some subsets of the logic space}

In this section we fix a countable continuous signature 
$$
L=\{ d,R_1 ,...,R_k ,..., F_1 ,..., F_l ,...\} , 
$$ 
a Polish space $({\bf Y},d)$.  a countable dense subset $S_{{\bf Y}}$ of ${\bf Y}$ 
and study subsets of  ${\bf Y}_L$ which are invariant with respect to 
isometries stabilising $S_{{\bf Y}}$ setwise. 
Viewing the logic space ${\bf Y}_L$ as a Polish space 
one can consider Borel/algorithmic complexity of 
some natural subsets of ${\bf Y}_L$ of this kind. 
This approach differs from the one of Section 2.4. 
It corresponds to considering a structure on ${\bf Y}$ 
(say $M$) together with its {\bf presentation over} $S_{{\bf Y}}$, 
i.e. the set 
$$ 
Diag (M,S_{{\bf Y}})=\{ (\phi ,q): M\models \phi <q, \mbox{ where } q\in [0,1]\cap \mathbb{Q} 
\mbox{ and } \phi \mbox{ is a continuous } 
$$ 
$$ 
\mbox{ sentence with parameters from }S_{{\bf Y}} \}.     
$$ 
It is natural in the cases of examples {\bf (A)} - {\bf (C)} of Section 2.3 and 
the corresponding computable presentations as in Section 3. 
Moreover it corresponds to the approach of computable model theory. 

We will concentrate on separable categoricity.

A theory $T$ is {\bf separably categorical} if any 
two separable models of $T$ are isomorphic. 
A useful reformulation of this notion is given in 
Theorem \ref{SCHm}. 
Since we will only use this theorem below 
all necessary facts concerning separable categoricity 
(together with the proof of Theorem \ref{SCHm}) 
are given in Appendix.

\subsection{Separable categoricity}

We preserve all the assumptions of Section 1 on 
the space $({\bf Y},d)$. 
For simplicity we assume that all $L$-symbols are 
of continuity modulus $\mathsf{id}$.  
Simplifying notation we put $S =S_{{\bf Y}}$.  
We reformulate separable categoricity as follows.

\begin{thm} \label{SCHm} 
Let $M$ be a non-compact, separable, continuous, 
metric structure on $({\bf Y},d)$.
The structure $M$ is separably categorical if and only 
if  for any $n$ and $\varepsilon$ there are finitely many  
conditions $\phi_i (\bar{x})\le \delta_i$, $i\in I$,  
so that any $n$-tuple of $M$ satisfies one of these 
conditions and the following property holds: 
\begin{quote} 
for any  
$i\in I$, and any $a_1 ,...,a_n \in M$ 
realizing  $\phi_i (\bar{x})\le \delta_i$ and  
any finite set of formulas $\Delta (x_1 ,...,x_n , x_{n+1})$ 
realized in $M$ and containing  $\phi_i (\bar{x})\le \delta_i$, 
there is a tuple $b_1 ,...,b_n ,b_{n+1}$ 
realizing $\Delta$ such that $max_{i\le n} d(a_i ,b_i ) <\varepsilon$.   
\end{quote} 
\end{thm} 

We now introduce a class of structures 
which is justified  by its formulation. 
If we assume that all parameters appearing in it can be taken from $S$ we 
arrive at the following definition. 

\begin{definicja} \label{DefSC}
 Let ${\cal SC}_S$  be the set of all $L$-structures $M$ on ${\bf Y}$ 
with the following condition: \\ 
{\em for every $n$ and rational $\varepsilon$ there is
a finite set $F$ of tuples $\bar{a}_i$ from $S$ together 
with conditions $\phi_i (\bar{x})\le \delta_i$  ($i\in I$ and all $\delta_i$ are rational) 
with $\phi^{M}_i (\bar{a}_i )\le \delta_i$, $i\in I$, 
and the following properties     
\begin{itemize} 
\item 
any $n$-tuple $\bar{a}$ from $S$ satisfies in $M$ one of these  
$\phi_i (\bar{x})\le \delta_i$ 

\item
when $\phi^{M}_i (\bar{a})\le \delta_i$ and $\bar{c}$ is an 
$(n+1)$-tuple from $S$ with $c_1 ,...,c_n$ satisfying 
$\phi_i (\bar{x})\le \delta_i$ in $M$, 
\begin{quote} 
for any finite set $\Delta$ of $L$-formulas 
$\phi (\bar{y})$, $|\bar{y}| =n+1$ with $\phi^{M} (\bar{c})=0$ 
there is an $(n+1)$-tuple  
$\bar{b}=(b_{1},\ldots b_{n+1})\in S$ so that \\
$\mathsf{max}_{j\le n} (d(a_j ,b_{j} ))\le \varepsilon$ and 
$ \phi^{M}(\bar{b})=0$ 
for all formulas $\phi \in \Delta$. 
\end{quote}  
\end{itemize} }
\end{definicja}

Note that Theorem \ref{SCHm} 
implies that if $M$ is a separably categorical structure on ${\bf Y}$, 
there is a dense set $S'\subseteq {\bf Y}$ so that $M$ belongs to 
the corresponding set of $L$-structures ${\cal SC}_{S'}$.   
To see this we just extend $S$ to some countable $S'$ 
which satisfies the property of  Definition \ref{DefSC}  
in which we additionally require that $a_1 ,...,a_n \in S'$. 
Thus the following statement becomes interesting.

\begin{prop} \label{SCBorel} 
The subset ${\cal SC}_S \subset {\bf Y}_L$  is  $\mathsf{Iso} (S)$-invariant 
and Borel. 
\end{prop}  

{\em Proof.} 
It is clear that  $\mathcal{SC}_S$ is $\mathsf{Iso} (S)$-invariant. 
To see that $\mathcal{SC}_S$ is a Borel subsets of ${\bf Y}_L$  
it suffices to note that given rational $\varepsilon >0$, 
finitely many formulas $\phi_i (\bar{x})$, $i\in I$,  
with $|\bar{x}|=n+1$, and an $n$-tuple  $\bar{a}$  from $S$  
the set of $L$-structures $M$ on ${\bf Y}$ with the property  that 
\begin{quote} 
there is an $(n+1)$-tuple $\bar{b}\in S$ so that 
$\mathsf{max}_{j\le n} (d(a_j ,b_j ))\le \varepsilon$ and 
$\phi^M_i (\bar{b})=0$ for all $i\in I$,  
\end{quote}  
is a Borel subset of ${\bf Y}_L$. 
The latter follows from Lemma \ref{EsLo}, which 
in particular says that any set of $L$-structures of the form 
$$
\{ M: M \models  \mathsf{max} (\mathsf{max}_{j\le n}(d(a_j ,b_j)\dot{-} \varepsilon ), \mathsf{max}_{i\in I} (\phi_i (\bar{b}) ))=0 \} 
$$
is a Borel subset of ${\bf Y}_L$. 
$\Box$ 

\bigskip 

The proof above demonstrates that $\mathcal{SC}_S$ is of Borel level $\omega$. 

In Section 4.3 we will discuss the conjecture that any separably categorical continuous 
$L$-structure on ${\bf Y}$ is homeomorphic to a structure from 
$\mathcal{SC}_S$.  
 
Since $\mathcal{SC}_S$ is a subset of the standard space ${\bf Y}_L$ we do need to specify 
grey bases $\mathcal{R}$, $\mathcal{U}$ as in Section 3.2. 
To be definite one can generate $\mathcal{R}$ by grey stabilizers and $\mathcal{U}$ 
by atomic formulas. 
This issue becomes important in the next section when we consider computable members of 
$\mathcal{SC}_S$. 

\subsection{Computable members} 

The following proposition is an effective version 
of Proposition \ref{SCBorel} in the case {\bf (A)} of Section 2.3. 
We work in the effective presentation given in Section 3.2 ({\bf D}). 

\begin{prop}  
Let  ${\cal SC}_{\mathbb{Q}\mathfrak{U}}$  be the 
$\mathsf{Iso} (\mathbb{Q}\mathfrak{U})$-invariant Borel subset of $\mathfrak{U}_L$ 
defined as in Section 4.1. 

Then the subset of indices of computable structures from 
${\cal SC}_{\mathbb{Q}\mathfrak{U}}$ 
is hyperarithmetical. 
\end{prop} 

{\em Proof.} 
Under the framework of Sections 3.1 and 3.2 ({\bf D}) the following 
statement holds. 
\begin{quote} 
The set of all pairs $(i,j)$ where $j$ is an index of a cone 
from $(\mathcal{B}_0 )_{\mathbb{Q}}$ and $i$ is an index 
of a computable structure from this cone, is hyperarithmetical 
of level $\omega$. 
\end{quote}   
This is an effective version of Proposition \ref{EsLo}. 
It follows from Lemma \ref{21}  by standard arguments. 
Note that as we have shown in Section 3.2 ({\bf D}) 
all assumptions of Lemma \ref{21} are satisfied under the circumstances of 
our proposition. 

It remains to verify that the definition of ${\cal SC}_{\mathbb{Q}\mathfrak{U}}$ 
defines a hyperarithmetical subset of indices of computable structures. 
This is straghtforward  
(similar to the proof of Proposition \ref{SCBorel}).  
$\Box$

\subsection{Countable dense homogeneity} 

We conjecture that when $({\bf Y}, d)$ is as in the cases {\bf (A)} - {\bf (C)} of 
Section 2.3 and the dense subset $S$ is chosen as the corresponding 
approximating substructure, then any separably categorical structure from 
${\bf Y}_L$ is homeomorphic to an element of $\mathcal{SC}_S$. 
In this section we connect it with  countable dense homogeneity.    

A separable space $X$ is {\bf countable dense homogeneous} (CDH) 
if given any two countable dense subsets $D$ and $E$ of $X$ there 
is a homeomorphism $f:X \rightarrow X$ such that 
$f(D) = E$ (see \cite{HrusMill}). 
It is known that the unbounded Urysohn space and the spaces $l^2$ 
are CDH (and they are homeomorphic). 

In \cite{dijkstra} J. Dijkstra introduced Lipschitz CDH as follows. 

\begin{definicja} 
A metric space $(X,d)$ is called Lipschitz countable dense homogeneous if 
given $\varepsilon$ and $A$, $B$ countable dense subsets of $X$ 
there is a homeomorphism $h : X\rightarrow X$ such that $f(A) = B$ and 
$$ 
1-\varepsilon < \frac{d(h(x), h(y))}{d(x,y))} < 1+\varepsilon \mbox{ for all } x,y\in X .
$$ 
\end{definicja} 

He has proved in \cite{dijkstra} that every separable Banach space 
is Lipschitz CDH.  
Moreover it is shown in \cite{KSh} that the unbounded 
Urysohn space is also Lipschitz CDH. 

Regarding the property $\mathcal{SC}_S$ this notion seems very helpful. 
Indeed, as we have already noted  Theorem \ref{SCHm} 
implies that if $M$ is a separably categorical structure on ${\bf Y}$, 
there is a dense set $S'\subseteq {\bf Y}$ so that $M$ belongs to 
the corresponding Borel set of $L$-structures ${\cal SC}_{S'}$.   
To support the conjecture of this subsection we need a homeomorphism 
which takes $S'$ onto $S$  and takes $M$ into $\mathcal{SC}_S$. 
However the latter condition is not easy to control.

\subsection{Complexity of sets of approximately ultrahomogeneous structures} 

If in the definition of the class $\mathcal{SC}_S$ we restrict ourselves by 
only quantifier free formulas we arrive at a definition 
of a subset of $\mathcal{SC}_S$ which we denote by $\mathcal{SCU}_S$. 
The approach of Sections 4.1 and 4.2 works in this case too. 
As in Section 4.3 one can conjecture that when 
$({\bf Y}, d)$ is as in the cases {\bf (A)} - {\bf (C)} of 
Section 2.3 and the dense subset $S$ is chosen as the corresponding 
approximating substructure, then any separably categorical ultrahomogeneous 
structure from ${\bf Y}_L$ is homeomorphic to an element of $\mathcal{SCU}_S$. 

It is worth noting here that 
since any Polish group can be realized as 
the automorphism group of an approximately ultrahomogeneous 
structure it makes sense  
to study Polish groups by description of the corresponding classes of approximately 
ultrahomogeneous structures and to study the complexity 
of these classes. 
For example we do not know if the class of approximately 
ultrahomogeneous $L$-structures on ${\bf Y}$ is 
a Borel subset of ${\bf Y}_L$.

\subsection{CLI Polish groups} 

In this section we give a different example of complexity of a subclass of ${\bf Y}_L$ 
for a Polish space  $({\bf Y}, d)$.  
It somehow corresponds to the result 
of M.Malicki \cite{mal} that the set of all Polish groups admitting 
compatible complete left-invariant metrics is coanalytic 
non-Borel as a subset of a standard Borel space of Polish groups. 

The following statement is Lemma 9.1 of \cite{BNT}. 

\begin{quote} 
Let $G$ be the automorphism group of a continuous $L$-structure $M$ 
on the space $({\bf Y}, d)$.  
Then the group $G$ admits a compatible complete left-invariant metric (i.e. $\mathsf{Aut}(M)$ is CLI) 
if and only if  each $L_{\omega_1 \omega}$-elementary embedding of $M$ into itself is  surjective. 
\end{quote} 

What is the complexity of the class of structures from this proposition? 
We have some remark concerning this question. 

\begin{prop} \label{corclim} 
Let $\mathcal{L}$ be a countable fragment of $L_{\omega_1 \omega}$. 
The subset of ${\bf Y}_L$ consisting of  structures $M$ 
admitting  proper  $\mathcal{L}$-elementary embeddings into itself is analytic. 
$\mathsf{Iso} (S)$-invariant and coanalytic. 
\end{prop} 

{\em Proof.}  
Consider the extension of $L$ by a unary function $f$.  
All expansions of $L$-structures satisfying the property 
that $f$ is an isometry which preserves the values $\mathcal{L}$-formulas, form a Borel 
subset of the (Polish) space of all $L\cup \{ f\}$-structures on ${\bf Y}$. 

If $s\in S=S_{{\bf Y}}$ and $\varepsilon \in {\bf Q}\cap [0,1]$ then 
the condition that $f(S)$ does not intersect the 
$\varepsilon$-ball of $s$ is open. 
Thus the set of $L\cup \{ f\}$-structures with a proper 
embedding $f$ into itself, is Borel. 
The rest is easy. 
$\Box$

\section{Appendix. Proof of Theorem \ref{SCHm} } 

We need the following definition. 

\begin{definicja} 
Let $A\subseteq M$. 
A predicate $P:M^n \rightarrow [0,1]$ is definable in 
$M$ over $A$ if there is a sequence $(\phi_k (x) :k\ge 1 )$ 
of $L(A)$-formulas such that predicates interpreting 
$\phi_k (x)$ in $M$ converge to $P(x)$ uniformly in $M^n$. 
\end{definicja} 

Let $p(\bar{x})$ be a type of a theory $T$. 
It is called {\bf principal} if for every model $M\models T$, the 
predicate $\mathsf{dist}(\bar{x},p(M))$ is definable over $\emptyset$. 

By Theorem 12.10 of \cite{BYBHU} a complete theory $T$ 
is separably categorical if and only if 
for each $n>0$, 
every $n$-type $p$ is principal. 
Another property equivalent to separable categoricity states that 
for each $n>0$, the metric space $(S_n (T),d)$ is compact.  
In particular for every $n$ and every $\varepsilon$ there is 
a finite family of principal $n$-types $p_1 ,...,p_m$ so that 
their $\varepsilon$-neighbourhoods cover $S_n(T)$. 

In the classical first order logic a countable 
structure $M$ is $\omega$-categorical if and only 
if $\mathsf{Aut}(M)$ is an {\bf oligomorphic} permutation group, 
i.e. for every $n$, $\mathsf{Aut}(M)$ has finitely many orbits 
on $M^n$. 
In continuous logic we have the following modification.  

\begin{definicja} 
An isometric action of a group $G$ on a metric space 
$({\bf X},d)$ is said to be {\bf approximately oligomorphic} 
if for every $n\ge 1$ and $\varepsilon >0$ there is 
a finite set $F\subset {\bf X}^n$ such that 
$$ 
G\cdot F = \{ g\bar{x} : g\in G \mbox{ and } \bar{x}\in F\}
$$
is $\varepsilon$-dense in $({\bf X}^n,d)$. 
\end{definicja} 

Assuming that $G$ is the automorphism group of a non-compact 
separable continuous metric structure $M$, $G$ is approximately 
oligomorphic if and only if  the structure $M$ is separably 
categorical (C. Ward Henson, see Theorem 4.25 in \cite{scho}). 
It is also known that separably categorical structures are 
{\bf approximately homogeneous} in the following sense: 
if $n$-tuples $\bar{a}$ and $\bar{c}$ have the same types 
(i.e. the same values $\phi (\bar{a})=\phi (\bar{b})$ for 
all $L$-formulas $\phi$) then for every $c_{n+1}$ and 
$\varepsilon >0$ there is an tuple $b_1 ,...,b_n ,b_{n+1}$ 
of the same type with $\bar{c},c_{n+1}$, so that 
$d(a_i, b_i )\le \varepsilon$ for $i\le n$.  
In fact for any $n$-tuples $\bar{a}$ and $\bar{b}$ 
there is an automorphism $\alpha$ of $M$ such that 
$$
d(\alpha (\bar{c}),\bar{a})\le d(tp(\bar{a}),tp(\bar{c})) +\varepsilon . 
$$  
(i.e $M$ is {\bf strongly $\omega$-near-homogeneous} 
in the sense of Corollary 12.11 of \cite{BYBHU}).

To prove Theorem  \ref{SCHm} we start with the following observation. 

\begin{lem} \label{SCH} 
Let $M$ be a non-compact, separable, continuous, 
metric structure on $({\bf Y},d)$.
The structure $M$ is separably categorical if and only if 
for any $n$ and $\varepsilon$ there are finitely many  
conditions $\phi_i (\bar{x})\le \delta_i$, $i\in I$,  
so that any $n$-tuple of $M$ satisfies one of these 
conditions and the following property holds: 
\begin{quote} 
for any  $i\in I$, any $a_1 ,...,a_n \in M$ 
realising  $\phi_i (\bar{x})\le \delta_i$ and any type 
$p(x_1 ,...,x_n , x_{n+1})$ realized in $M$ and containing  
$\phi_i (\bar{x})\le \delta_i$, there is a tuple 
$b_1 ,...,b_n ,b_{n+1}$ realizing $p$ such that 
$max_{i\le n} d(a_i ,b_i ) <\varepsilon$.   
\end{quote} 
\end{lem} 

{\em Proof.} 
By Theorem 12.10 of \cite{BYBHU} a complete theory $T$ is separably categorical 
if and only if for each $n>0$, every $n$-type is principal. 
An equivalent condition states that for each $n>0$, the metric space 
$(S_n (T),d)$ is compact.  
In particular for every $n$ and every $\varepsilon$ there is 
a finite family of principal $n$-types $p_1 ,...,p_m$ so that 
their $\varepsilon/2$-neighbourhoods cover $S_n(T)$. 

Thus when $M$ is separably categorical, given $n$ and $\varepsilon$, 
we find appropriate $p_i$, $i\in I$, define    
$P_i (\bar{x})=\mathsf{dist}(\bar{x}, p_i (M))$, 
the corresponding definable predicates and  
$n$-conditions $\phi_i (\bar{x})\le \delta_i$ describing the corresponding 
$\varepsilon /2$-neighbourhoods of $p_i$.  
The rest follows by strong $\omega$-near-homogeneity. 

To see the converse assume that $M$ satisfies the property from the 
formulation. 
To see that $G=Aut (M)$ is approximately oligomorphic 
take any $n$ and $\varepsilon$ and find finitely many  
conditions $\phi_i (\bar{x})\le \delta_i$, $i\in I$,  
satisfying the property from the formulation for $n$ and $\varepsilon /4$. 
Choose $\bar{a}_i$ with $\phi_i (\bar{a}_i )\le \delta_i$ and 
let $F =\{ \bar{a}_i : i\in I\}$.  
To see that $G\cdot F$ is $\varepsilon$-dense 
we only need to show that if $\bar{a}$ satisies 
$\phi_i (\bar{x})\le \delta_i$, then there is an automorphism 
which takes $\bar{a}$ to the $\varepsilon$-neighbourhood 
of $\bar{a}_i$. 
This is verified by "back-and-forth" as follows. 
Let $(\varepsilon_k )$ be an infinite sequence of positive 
real numbers whose sum is less than $\varepsilon /4$. 
At every step $l$ (assuming that $l\ge n$) we build 
a finite elementary map $\alpha_l$ and $l$-tuples 
$\bar{c}_l$ and $\bar{d}_l$ so that  
\begin{itemize} 
\item $\bar{c}_n =\bar{a}$ and $\bar{d}_n =\bar{a}_i$;  
\item for $l>n$, $\alpha_l$ takes $\bar{c}_l$ to $\bar{d}_l$
\item for $l>n+1$, the first $l-1$ coordinates of $\bar{c}_l$ 
(resp. $\bar{d}_l$) are at distance less than  $\varepsilon_l$ 
away from the corresponding coordinates of $\bar{c}_{l-1}$ 
(resp. $\bar{d}_{l-1}$); 
\item the sets $\bigcup\{ \bar{c}_l :l\in \omega\}$ and 
$\bigcup\{ \bar{d}_l: l\in\omega\}$ are dense in $M$. 
\end{itemize} 
In fact we additionally arrange that for even $l$, 
$\bar{c}_{l+1}$ extends $\bar{c}_l$ and for odd $l$ 
$\bar{d}_{l+1}$ extends $\bar{d}_l$. 
In particular the type of  $\bar{c}_{l+1}$ 
always extends the type of $\bar{c}_l$.   
At the $(n+1)$-th step we find finitely many conditions 
$\phi'_j (\bar{x})\le \delta'_j$, $j\in J$,  
so that any $(n+1)$-tuple of $M$ satisfies one of these 
conditions and for any  $j\in J$, any $a'_1 ,...,a'_{n+1} \in M$ 
realising  $\phi'_j (\bar{x})\le \delta'_j$ and any type 
$p(x_1 ,...,x_{n+1} , x_{n+2})$ realized in $M$ and 
containing  $\phi'_j (\bar{x})\le \delta_j$, there is 
a tuple $b_1 ,...,b_{n+1} ,b_{n+2}$ realizing $p$ such 
that $\mathsf{max}_{t\le n+1} d(a'_t ,b_t ) <\varepsilon_{n+1}$. 
Now by the choice of $i$ for any extension 
of $\bar{a}=\bar{c}_{n}$ to an $(n+1)$-tuple $\bar{c}_{n+1}$ 
we can find a tuple $\bar{d}_{n+1}$ realizing $tp(\bar{c}_{n+1})$  
so that the first $n$ coordinates of $\bar{d}_{n+1}$ are 
at distance less than $\varepsilon /4$ away from 
the corresponding coordinates of $\bar{d}_{n}=\bar{a}_i$.  
If $n$ is even we choose such $\bar{c}_{n+1}$ and 
$\bar{d}_{n+1}$; if $n$ is odd we replace the roles 
of $\bar{c}_{n+1}$ and $\bar{d}_{n+1}$. 
For the next step we fix the condition 
$\phi'_j (\bar{x})\le \delta'_j$ satisfied 
by $\bar{c}_{n+1}$ and $\bar{d}_{n+1}$. 

The $(l+1)$-th step is as follows. 
Assume that $l$ is even (the odd case is symmetric). 
Extend $\bar{c}_l$ to an appropriate $\bar{c}_{l+1}$ 
(aiming to density of $\bigcup \{ \bar{c}_l :l\in \omega\}$). 
There are finitely many conditions $\phi''_k (\bar{x})\le \delta''_k$, 
$k\in K$, so that any $(l+1)$-tuple of $M$ satisfies one of these 
conditions and for any  $k\in K$, any $a'_1 ,...,a'_{l+1} \in M$ 
realising  $\phi''_k (\bar{x})\le \delta''_k$ and any type 
$p(x_1 ,...,x_{l+1} , x_{l+2})$ realized in $M$ and 
containing  $\phi''_k (\bar{x})\le \delta''_k$, there is 
a tuple $b_1 ,...,b_{l+1} ,b_{l+2}$ realizing $p$ such that 
$max_{t\le l+1} d(a'_t ,b_t ) <\varepsilon_{l+1}$.   
We find the condition satisfied by $\bar{c}_{l+1}$ 
and a tuple $\bar{d}_{l+1}$ realizing $tp(\bar{c}_{l+1})$  
so that the first $l$ coordinates of $\bar{d}_{l+1}$ are at distance 
less than $\varepsilon_{l}$ away from the corresponding 
coordinates of $\bar{d}_{l}$. 

As a result for every $k$ we obtain Cauchy sequences of 
$k$-restrictions of $\bar{c}_l$-s and $\bar{d}_l$-s. 
For $k=n$ their limits are not distant from $\bar{a}$ 
and $\bar{a}_i$ more than $\varepsilon/2$. 
Moreover the limits $lim \{\bar{c}_l \}$ and $lim \{\bar{d}_l \}$ 
are dense subsets of ${\bf Y}$ and realize the same type. 
This defines the required automorphism of $M$. 
$\Box$

\bigskip 

{\em Proof of Theorem \ref{SCHm}.} 
It suffices to show that the condition of the 
formulation implies the corresponding condition of Lemma \ref{SCH}.  
Given $n$ and $\varepsilon$ take the family 
$\phi_i (\bar{x})\le \delta_i$, $i\in I$, satisfying 
the condition of the theorem for $n$ and $\varepsilon /2$. 
Let $p(\bar{x},x_{n+1})$ be a type with  
$\phi_i (\bar{x})\le \delta_i$ and $a_1 ,...,a_n$ 
be as in the formulation. 

Let $(\varepsilon_k )$ be an infinite sequence of positive 
real numbers whose sum is less than $\varepsilon /2$. 
Now apply the condition of the formulation of the theorem 
to $n+1$ and $\varepsilon_1 /2$ and find an appropriate 
finite family of inequalities such that one of them, 
say  $\psi (\bar{x},x_{n+1}) \le \tau$, belongs to $p$ 
and for any $c_1 ,...,c_n ,c_{n+1} \in M$ 
realising  $\psi (\bar{x},x_{n+1})\le \tau$, 
and any finite subset $\Delta \subset p$ containing  
$\psi (\bar{x},x_{n+1})\le \tau$
there is a tuple $c'_1 ,...,c'_n ,c'_{n+1}$ realizing  $\Delta$,  
such that $\mathsf{max}_{i\le n+1} d(c_i ,c'_i ) <\varepsilon_1 /2$. 
Then let    $b^1_1 ,...,b^1_n ,b^1_{n+1}$ be a tuple realizing   
$\phi_i (\bar{x})\le \delta_i$ and 
$\psi (\bar{x},x_{n+1})\le \tau$ such that 
$\mathsf{max}_{i\le n} d(a_i ,b^1_i ) <\varepsilon /2$. 

For $n+1$ and $\varepsilon_2 /2$ find an appropriate condition 
$\psi' (\bar{x},x_{n+1}) \le \tau'$ from $p$ 
so that  any $c_1 ,...,c_n ,c_{n+1} \in M$ realizing  
$\psi' (\bar{x},x_{n+1})\le \tau'$, 
and any finite subset $\Delta \subset p$ containing  
$\psi' (\bar{x},x_{n+1})\le \tau'$
there is a tuple $c'_1 ,...,c'_n ,c'_{n+1}$ realizing  $\Delta$,  
such that $\mathsf{max}_{i\le n+1} d(c_i ,c'_i ) <\varepsilon_2 /2$. 
Let    $b^2_1 ,...,b^2_n ,b^2_{n+1}$ be a tuple 
realising  $\phi_i (\bar{x}) <\delta_i$,  
$\psi (\bar{x},x_{n+1})\le \tau$ and 
$\psi' (\bar{x},x_{n+1})\le \tau'$ such that 
$\mathsf{max}_{i\le n+1} d(b^1_i ,b^2_i ) <\varepsilon_1 /2$. 
Note that  $\mathsf{max}_{i\le n} d(a_i ,b^2_i ) <\varepsilon /2 +\varepsilon_1 /2$. 

Continuing this procedure we obtain a Cauchy sequence of 
$(n+1)$-tuples so that its limit satisfies $p$ 
and is not distant from $\bar{a}$ more than $\varepsilon$. 
$\Box$

\bigskip

Institute of Mathematics, Silesian Univesity of Technology, \parskip0pt

ul. Kaszubska 23, 44-101 Gliwice, Poland \parskip0pt

E-mail: Aleksander.Iwanow@polsl.pl

\bigskip 

Institute of Mathematics and Computer Science, \parskip0pt 

University of Opole, ul.Oleska 48, 45 - 052 Opole, Poland \parskip0pt 

E-mail: bmajcher@uni.opole.pl   
\end{document}